\documentclass[12pt]{article}

\usepackage[latin1]{inputenc}
\usepackage{bbm}
\usepackage{stmaryrd}
\usepackage{latexsym}
\usepackage{amsfonts}
\usepackage{amsmath,amssymb,amscd}
\usepackage{yfonts}
\usepackage{dsfont}
\usepackage{mathabx}
\usepackage[dvips]{graphicx}
\usepackage{epsfig}
\usepackage{psfrag}
\usepackage[font={small}]{caption}
\usepackage{color}
\usepackage{float}
\usepackage{upgreek}
\usepackage{array}

\DeclareFontFamily{U}{txsyc}{}
\DeclareFontShape{U}{txsyc}{m}{n}{
   <-> txsyc%
}{}
\DeclareFontShape{U}{txsyc}{bx}{n}{
   <-> txbsyc%
}{}
\DeclareFontShape{U}{txsyc}{l}{n}{<->ssub * txsyc/m/n}{}
\DeclareFontShape{U}{txsyc}{b}{n}{<->ssub * txsyc/bx/n}{}
\DeclareSymbolFont{symbolsC}{U}{txsyc}{m}{n}
\SetSymbolFont{symbolsC}{bold}{U}{txsyc}{bx}{n}
\DeclareFontSubstitution{U}{txsyc}{m}{n}
\DeclareMathSymbol{\df}{\mathrel}{symbolsC}{"42}
\DeclareMathSymbol{\fd}{\mathrel}{symbolsC}{"43}
\DeclareMathSymbol{\lJoin}{\mathrel}{symbolsC}{"58}
\DeclareMathSymbol{\rJoin}{\mathrel}{symbolsC}{"59}

\newcommand{\G}{G}
\newcommand{\V}{V}
\newcommand{\E}{E}
\newcommand{\U}{\mathcal{U}}

\usepackage{latexsym}
\newcommand{\cA}{{\cal A}}

\newcommand{\cC}{{\cal C}}

\newcommand{\cE}{{\cal E}}
\newcommand{\cF}{{\cal F}}

\newcommand{\cH}{{\cal H}}

\newcommand{\cK}{{\cal K}}
\newcommand{\cL}{{\cal L}}
\newcommand{\cN}{{\cal N}}

\newcommand{\cP}{{\cal P}}

\newcommand{\CC}{\mathbb{C}}

\newcommand{\EE}{\mathbb{E}}

\newcommand{\LL}{\mathbb{L}}
\newcommand{\NN}{\mathbb{N}}
\newcommand{\PP}{\mathbb{P}}

\newcommand{\RR}{\mathbb{R}}

\newcommand{\ZZ}{\mathbb{Z}}

\newcommand{\iy}{\infty}
\newcommand{\lt}{\left}
\newcommand{\me}{\medskip}

\newcommand{\pa}{\partial}
\newcommand{\ri}{\rightarrow}
\newcommand{\rt}{\right}
\newcommand{\sm}{\smallskip}

\newcommand{\wi}{\widetilde}
\newcommand{\wit}{\widehat}
\newcommand{\card}{\mathrm{card}}

\newcommand{\ex}{\exists\ }
\newcommand{\fo}{\forall\ }

\newcommand{\lve}{\lt\vert}
\newcommand{\lVe}{\lt\Vert}

\newcommand{\rve}{\rt\vert}
\newcommand{\rVe}{\rt\Vert}

\newcommand{\st}{\,:\,}

\newcommand{\trace}{\mathrm{tr}}

\newcommand{\un}{\mathds{1}}

\newcommand{\vvv}{\vert\!\vert\!\vert}
\newcommand{\bq}{\begin{eqnarray*}}
\newcommand{\bqn}[1]{\begin{eqnarray}\label{#1}}
\newcommand{\eq}{\end{eqnarray*}}
\newcommand{\eqn}{\end{eqnarray}}
\newcommand{\wwtbp}{\par\hfill $\blacksquare$\par\me\noindent}
\newcommand{\thistitlepagestyle}{}
\newcommand{\lin}{\llbracket}
\newcommand{\rin}{\rrbracket}

\newcommand{\ttsim}{\raise.17ex\hbox{$\scriptstyle\mathtt{\sim}$}}

\newtheorem{pro}{Proposition}
\newtheorem{cor}[pro]{Corollary}
\newtheorem{lem}[pro]{Lemma}
\newtheorem{theo}[pro]{Theorem}
\renewcommand{\thepro}{\arabic{pro}}

\newenvironment{rem}
{\par\me\refstepcounter{pro}\noindent{\bf Remark \thepro\ }}
{\par\hfill $\square$\par\sm\noindent}

\newcommand{\proof}{\par\me\noindent\textbf{Proof}\par\sm\noindent}
\newcommand{\prooff}[1]{\par\me\noindent\textbf{#1}\par\sm\noindent}

\newcommand{\fX}{\mathfrak{X}}
\newcommand{\fF}{\mathfrak{F}}

\title{On the fastest finite Markov processes\footnote{Work supported by the Indo-French Centre for Applied Mathematics and by the ANR-12-BS01-0019.}}
\author{Vivek Borkar${}^\dagger$ and Laurent Miclo${}^\ddagger$
}
\date{\vbox{\copy0
\vskip5mm
\copy1
}
}

\begin{document}

\setbox0=\vbox{
\large
\begin{center}
 ${}^\dagger$Department of Electrical Engineering\\
 Indian Institute of Technology Bombay, India
\end{center}
}
 \setbox1=\vbox{
 \large
 \begin{center}
 ${}^\ddagger$Institut de Mathématiques de Toulouse, UMR 5219\\
Université de Toulouse and CNRS, France\\
 \end{center}
}
\setbox5=\vbox{
\hbox{${}^\dagger$borkar@ee.iitb.ac.in\\}
\vskip1mm
\hbox{Department of Electrical Engineering\\}
\hbox{Indian Institute of Technology Bombay\\}
\hbox{Mumbai 400076, India\\}
}

\setbox6=\vbox{
\hbox{${}^\ddagger$miclo@math.univ-toulouse.fr\\}
\vskip1mm
\hbox{Institut de Mathématiques de Toulouse\\}
\hbox{Université Paul Sabatier\\}
\hbox{118, route de Narbonne\\}
\hbox{31062 Toulouse Cedex 9, France\\}
}

\maketitle
\thistitlepagestyle
\abstract{Consider a finite irreducible Markov chain with invariant probability $\pi$.
Define its inverse communication speed as the expectation to go from $x$ to $y$, when $x,y$
are sampled independently according to $\pi$. In the discrete time setting and when $\pi$ is the uniform distribution $\upsilon$, Litvak and Ejov \cite{MR2542990}
have shown that the permutation matrices associated to Hamiltonian cycles are the fastest Markov chains.
Here we prove (A) that the above optimality is with respect to all processes compatible with a fixed graph of permitted transitions (assuming that it does contain a Hamiltonian cycle), not only the Markov chains, and, (B) that this result admits a natural extension in both discrete and continuous time when $\pi$ is close to $\upsilon$:
the fastest Markov chains/processes
are those moving successively on the points of a Hamiltonian cycle, with transition probabilities/jump rates dictated by $\pi$.
Nevertheless, the claim is no longer true when $\pi$ is significantly different from $\upsilon$.

}
\vfill\null
{\small
}\par

\newpage

\section{Introduction}

Given a finite oriented (strongly) connected graph $G=(V,E)$ and a positive probability measure $\pi$ on $V$,
it is natural to wonder what is the fastest Markov chain leaving $\pi$ invariant and whose permitted transitions are included in $E$.
This depends on the way the speed is measured. In this paper the goal is to minimize the expectation $\fF$ of the time needed to go from $x$ to $y$,
when $x$ and $y$ are independently sampled according to $\pi$.
Litvak and Ejov \cite{MR2542990} have shown that if $\pi$ is the \textbf{uniform distribution} $\upsilon$ and if $G$ contains a Hamiltonian cycle,
then the fastest Markov chains are exactly those following deterministically the succession of the states given by a Hamiltonian cycle when one exists
(the corresponding  quantity $\fF$ does not depend on the choice of the admissible Hamiltonian cycle).
Our objectives in this paper are: (A) to extend this result to the continuous time framework (under an appropriate renormalization of the jump rates), (B) to establish the above optimality over a larger class of processes,
and  to begin an investigation of the situation where $\pi$ is not the uniform distribution by showing, (C) that when $G$ contains a Hamiltonian cycle  and that $\pi$ is close to $\upsilon$,  the fastest Markov chains/processes are still those appropriately associated to Hamiltonian cycles, and (D) that this is no longer true when $\pi$ is `far away' from $\upsilon$.\\

The plan of the paper is as follows. The above results (A) and (B)  are proved in the next section via a dynamic programming approach,
which also provides an alternative proof of the discrete time result of Litvak and Ejov \cite{MR2542990}.
In Section 3, we decompose the generators leaving $\pi$ invariant  into convex sums of generators associated to (not necessarily Hamiltonian) cycles and we differentiate the expectations of hitting times with respect to the generators.
This is the basic tool for the proof of (C) (see Theorem \ref{theo2} in Section 3) in Section 4 , through small perturbations of the uniform probability measure. At the other extreme, large perturbations
lead to the proof of (D) (cf.\ Theorem \ref{theo3} in Section 3) at the end of the same section.
Section 5 contains some observations about the links between continuous time and discrete time.
In the appendix, we compute the fastest normalized birth and death generators leaving invariant any fixed positive probability measure $\pi$ on $\{0,1,2\}$.
The underlying graph
 is the segment graph of  length 2, i.e.\ the simplest example not containing a Hamiltonian cycle.\\

\section{The dynamic programming approach}

\subsection{Introduction}

The aim of this section is to show that the Hamiltonian cycles, when one exists, are the fastest in the sense we have defined among \textit{all} processes compatible with the given graph, not just the Markov chains. The proof uses dynamic programming. We first recall the eigentime identity in the next subsection and then establish the desired result for resp.\ discrete and continuous time in the subsections that follow.\\

\subsection{The eigentime identity}

We shall use the notation $\mathcal{L}(X)$ to denote the law of a random variable $X$ and $|A|$ for the cardinality of a finite set $A$. Consider a discrete time  Markov chain $(X_n)_{n\in\ZZ_+}$ on a finite state space $V$  with transition matrix $P = (p(i,j))_{i,j \in V}$. We assume it to be \textbf{irreducible}, i.e., for any $i,j \in V$, there exists a path $i_0 = i, i_1, \cdots, i_{n-1}, i_n = j$ such that $p(i_k, i_{k+1}) > 0$ , for $k\in\lin 0,n-1\rin\df \{0,1, ...,n-1\}$. Let $\pi \df (\pi(i))_{i\in V}$ denote its unique \textbf{stationary distribution}, which is its left eigenvector corresponding to the Perron-Frobenius eigenvalue $\theta_1 = 1$.
 In particular, if $\mathcal{L}(X_0) = \pi$, then for any $n\in\ZZ_+$, we have $\mathcal{L}(X_n) = \pi$. This justifies the term `stationary', its uniqueness being a well-known consequence of the irreducibility hypothesis. Denote by $\theta_2, \cdots, \theta_{\vert V\vert}$ the remaining eigenvalues of $P$. Also define the hitting times
\bqn{hitting}
\fo i\in V,\qquad \uptau_i &\df& \min\{n\in\ZZ_+\st  X_n = i\}.
\eqn
The \textbf{eigentime identity} states that
\begin{equation}
\forall \ i \in V,\qquad
\sum_j\pi(j)E_i[\uptau_j]\ =\ \sum_{m =2}^{\vert V\vert}\frac{1}{1 - \theta_m}. \label{eigen1}
\end{equation}
Here each eigenvalue is counted as many times as its (algebraic) multiplicity. For reversible chains, this is Proposition 3.13, p.\ 75, of Aldous and Fill \cite{aldous-fill-2014}. It was extended to the general case in Cui and Mao \cite{MR2728712}, see also \cite{MR3442821} for a simple proof and further extensions. The left hand side gives the mean hitting time of a target state picked randomly with distribution $\pi$. A minor modification is to consider instead the stopping times
$$\fo  i \in V,\qquad T_i\ \df\ \min\{n\in\NN\st X_n = i\},$$
where $\NN\df\ZZ_+\setminus\{0\}$, the set of positive integers.
Since $\pi(i)E_i[T_i] = 1$, we can replace (\ref{eigen1}) by
\begin{equation}
 \forall \ i \in V,\qquad \sum_j\pi(j)E_i[T_j]\ =\ 1 + \sum_{m =2}^{\vert V\vert}\frac{1}{1 - \theta_m}. \label{eigen2}
\end{equation}
This variant is essentially contained in Theorem 2.4 of Hunter \cite{Hunter}, which also gives the pre--history of the problem going back to Kemeny and Snell \cite{Kemeny}. An immediate  corollary is the `symmetrized' version
\begin{eqnarray}
\sum_{i,j}\pi(i)\pi(j)E_i[\uptau_j] &=& \sum_{m =2}^{\vert V\vert}\frac{1}{1 - \theta_m}, \label{aveigen1} \\
\sum_{i,j}\pi(i)\pi(j)E_i[T_j] &=& 1 + \sum_{m =2}^{\vert V\vert}\frac{1}{1 - \theta_m}. \label{aveigen2}
\end{eqnarray}

 Let  $\Pi $ be the rank one matrix  whose rows are all equal to $\pi$. From Theorem 2.4 of Hunter \cite{Hunter}, we have:\\

 (\ref{aveigen2}) equals tr$(I - P + \Pi)^{-1}$, where $I$ is the identity matrix. \hfill \textbf{$(*)$}\\

A \textbf{Hamiltonian cycle} $A$ of $V$ is an ordering $(a_0, a_1, ..., a_{N-1})$ of the elements of $V$, where $N\df\vert V\vert$.
We will make the convention that $a_N=a_0$,
as the indices should be seen as elements of $\ZZ_N\df\ZZ/(N\ZZ)$. More precisely, the cycles $(a_0, a_1, ..., a_{N-1})$ and $(a_k, a_{k+1}, ..., a_{k+N-1})$
should be identified as the same cycle, for all $k\in\ZZ_N$. This will be implicit in the sequel, even though for notational convenience, we will represent a cycle $A$ as $(a_0, a_1, ..., a_{N-1})$.
Consider an irreducible directed graph $\G = (\V, \E)$ where $\V, \E$ denote respectively its node and edge sets.
In the discrete time setting, such graphs will always be assumed to contain all the self-loops, i.e.\ $(i,i)\in E$ for any $i\in V$.
The Hamiltonian cycle $A\df (a_0, a_1, ..., a_{N-1})$ is said to be \textbf{admissible} for $G$ if $(a_k,a_{k+1})\in E$ for all $k\in \ZZ_n$.
It means that $G_A$ is a subgraph of $G$, where $G_A$ is the oriented graph on $V$ whose edges are the $(a_k,a_{k+1})$, for $k\in\ZZ_N$.
The set of all Hamiltonian cycles (respectively, admissible for $G$) is denoted $\cH$ (resp., $\cH(G)$) and the graph $G$ is said to be \textbf{Hamiltonian} if $\cH(G)\not =\emptyset$.
Obviously, we have $\cH=\cH(K_V)$, where $K_V$ is the complete oriented graph on $V$.

Consider the optimization problem of minimizing (\ref{eigen1})/(\ref{aveigen1}), or equivalently, (\ref{eigen2})/(\ref{aveigen2}), where $\pi$ is equal to $\upsilon$, over all irreducible $P$ compatible with the given graph $G$ (in the sense
that  for any $i\not =j\in V$, $P(i,j)>0 \Rightarrow (i,j)\in E$). Since these quantities will be infinite for reducible $P$, we might as well consider the problem of minimizing it over all stochastic matrices $P$ compatible with $\G$. Say that $P$ is Hamiltonian if there exists a Hamiltonian cycle $(a_1, \cdots, a_{\vert V\vert})$ such that $p(a_k, a_{k+1}) = 1 = p(a_{\vert V\vert}, a_1)$ for $1 \leq  k < \vert V\vert$. That is, the transitions deterministically trace a Hamiltonian cycle. Recall that a Hamiltonian cycle need not exist in general and the problem of determining whether one does is NP-hard (see, e.g., Garey and Johnson \cite{Garey}).  By Proposition 2.1 of Litvak and  Ejov \cite{MR2542990} in combination with \textbf{$(*)$} above, we have:\\

\begin{theo}\label{theo-1} When $\pi=\upsilon$, either of the quantities (\ref{eigen1}), (\ref{eigen2}) is minimized by a Hamiltonian $P$ if there exists one.
\end{theo}

In the next subsection, we give an alternative proof, inspired by the Held-Karp algorithm for scheduling problems \cite{Held}, which gives a  strengthening of this result and
has interesting implications for random search.
 Specifically, we improve on the cited result of  Litvak and  Ejov \cite{MR2542990}, insofar as the cost is shown to be minimized over \textit{all} $\G$-compatible random processes and not only the Markov chains, by tracing the Hamiltonian cycle, when one exists, deterministically.

\subsection{A dynamic programming solution}

As in Held and Karp \cite{Held}, a natural state space for the dynamic program is
$$V^* := \{(i,A) : i \in V, \ A \subset V\backslash\{i\}\}.$$
  With each $(i,A) \in V^*$, we associate an action space $U_i :=$ the set of probability vectors on the set
$$\V_i := \{j \in \V : (i,j) \in \E\} \subset \V$$
of successors of $i$ in $\G$. Note that this does not depend on $A$. Suppose $\V_i$ is enumerated as $(j_1, \cdots, j_{m_i})$. Given a `control' $q = (q(j_1), \cdots, q(j_{m_i}))$, the transition probability
\begin{equation}
\hat{p}((j, B) | (i, A), q) \label{trans}
\end{equation}
of going from $(i,A) \in V^*$ to $(j,B) \in V^*$ under control $q$ is zero if either $j \notin \V_i$ or $B \neq A\backslash\{j\}$. Otherwise it equals $q(j)$. Consider an $V^*$-valued controlled Markov chain $(X_n, Z_n)_{n\in\ZZ_+}$ governed by a control process $(q_n)_{n\in\ZZ_+}$ with $q_n \in U_{X_n}$, for all $n\in\ZZ_+$, evolving according to the above controlled transition probability function. That is, for any $n\in\ZZ_+$,
\begin{eqnarray*}
\lefteqn{P((X_{n+1}, Z_{n+1}) = (j, B) | (X_m, Z_m), q_m, m \leq n) } \\
&=& P((X_{n+1}, Z_{n+1}) = (j, B) | (X_n, Z_n), q_n) \\
&=& q_n(j)\delta_{B, A\backslash\{j\}},
\end{eqnarray*}
where $\delta_{\cdot , \cdot}$ denotes the Kronecker delta. Since we are allowed to choose any past dependent transition probability compatible with $\G$, this covers \textit{all} $V$-valued random processes that are compatible with $\G$, i.e., that make transitions only along the edges in $\E$. \\

Our objective is to minimize, for a prescribed initial state $i_0$\footnote{more generally, for a prescribed initial distribution} the quantity
\begin{equation}
E\left[\sum_j\uptau_j \ \Big| \ X_0= i_0, Z_0 = V\backslash\{i_0\}\right], \label{cost}
\end{equation}
which is proportional to (\ref{eigen1}) when  $\pi =\upsilon$, the uniform distribution (i.e., when $P$ is doubly stochastic). Note, however, that we do not require $(X_n)_{n\in\ZZ_+}$ to be even Markov. Let
$$\zeta := \min\{n \geq 0 : Z_n = \emptyset\}.$$
Then (\ref{cost}) can be equivalently written as
\begin{equation}
E\left[\sum_{m=0}^{\zeta}|Z_n| \ \Big| \ X_0 = i_0, Z_0 = V\backslash\{i_0\}\right]. \label{cost1}
\end{equation}
This allows us to apply the dynamic programming principle to the `value function' or `cost to go function'
$$V(i,A) := \inf E\left[\sum_{m=0}^{\zeta}|Z_n| \ \Big| \ X_0 = i, Z_0 = A\right],$$
where the infimum is over all admissible controls. Standard arguments yield the dynamic programming equation
\begin{eqnarray}
V(i,A) &=& \min_{q \in U_i}\left(|A| + \sum_{j \in \V_i}q(j)V(j, A\backslash\{j\})\right), \ A \neq \emptyset, \label{DP1} \\
V(\cdot, \emptyset) &\equiv& 0. \label{DP2}
\end{eqnarray}
Furthermore, the optimal control in state $(i,A)$ is any minimizer of the right hand side of (\ref{DP1}). Since the expression being minimized is affine in $q$, this minimum will be attained at a Dirac measure, implying that the optimal choice in state $(i,A)$ is to deterministically move to a certain $j \in \V_i$. In other words, the optimal trajectory is deterministic and perforce visits each node at least once, otherwise the cost would be infinite. Since at most one new node can be visited each time, the total cost is at least $\sum_{i=1}^{\vert V\vert-1}i = \frac{\vert V\vert(\vert V\vert-1)}{2}$, which equals the cost for tracing a Hamiltonian cycle if one exists.\\

A parallel treatment can be given for the cost
\begin{equation}
E\left[\sum_jT_j \ | \ X_0= i_0, Z_0 = V\right], \label{cost2}
\end{equation}
which can be equivalently written as
\begin{equation}
E\left[\sum_{m=0}^{\zeta}|Z_n| \ | \ X_0 = i_0, Z_0 = V\right]. \label{cost3}
\end{equation}
The minimum cost for this, again attained by tracing a Hamiltonian cycle deterministically, will be $\frac{s(s+1)}{2}$. \\

  We have proved:  \\

\begin{theo}\label{theo-2} Minimum of either the cost (\ref{cost}) or the cost (\ref{cost2}) over all $V$-valued random processes compatible with $\G$ is attained by tracing a Hamiltonian cycle when one exists.
\end{theo}

\ \\

This has interesting implications to some random search schemes. For example, consider the problem of searching for an N bit binary password given a device or `oracle' that can verify whether a  password is correct or not. Random search schemes for this problem have been proposed, involving Markov chains on the discrete $N$-cube $\{0, 1\}^N$, where any two strings differing in one position are deemed to be neighbors. This undirected graph can be rendered directed by replacing each undirected edge by two directed edges. A simple induction argument shows that it has a Hamiltonian cycle. Then the foregoing leads to the conclusion that no random search scheme can do better on average than simply listing the $N$-strings and checking them one by one.

\subsection{Continuous time problem}

 We now consider the continuous time counterparts of the foregoing.  Recall that a \textbf{Markov generator} on $V$ can be represented by a matrix $L\df(L(x,y))_{x,y\in V}$ whose off-diagonal entries are non-negative
and whose row sums all vanish. Corresponding Markov processes, defined through the corresponding martingale problems, will be denoted $X\df(X_t)_{t\geq 0}$.
The law of $X$ then only depends on the \textbf{initial distribution}, namely on the law $\cL(X_0)$ of $X_0$.
The Markov generator $L$ is said to be \textbf{compatible} with $G$, if we have
\bq
\fo x\not=y \in V,\qquad L(x,y)>0&\Rightarrow& (x,y)\in E.\eq
\par
The probability measure $\pi$, viewed as a row vector, is said to be \textbf{invariant} for the generator $L$, if $\pi L=0$.
Its probabilistic interpretation is that if initially $\cL(X_0)=\pi$, then for any $t\geq 0$,
$\cL(X_t)=\pi$, similarly to the discrete time case.
The generator $L$ is said to be \textbf{irreducible} if for any $x, y\in V$, there exists a path $x_0=x, x_1, ..., x_l=y$, with $l\in \ZZ_+$ the length of the path,
such that $L(x_k,x_{k+1})>0$ for all $k\in\lin 0,l-1\rin$.
In our finite setting, a Markov generator $L$ always admits an invariant probability measure, the irreducibility of $L$ ensures that it is unique.
The irreducible Markov generator $L$ is said to be \textbf{normalized}, if
\bqn{normalisation}
\sum_{x\in V} L(x) \pi(x)&=&1,\eqn
where $\mu$ is the invariant measure of $L$ and where $L(x)\df-L(x,x)=\sum_{y\not=x} L(x,y)$ for any $x\in V$.
It means that at its equilibrium $\mu$ (i.e.\ for the stationary $X$ starting with $\cL(X_0)=\mu$), the jump rate of $X$ is 1.
Denote by $\cL(G,\pi)$ the convex set of irreducible normalized Markov generators $L$ compatible with $G$  and admitting $\pi$ for invariant probability.
 To simplify notation, we will also write $\cL(\pi)\df \cL(K_V,\pi)$, when $G$ is the complete graph $K_V$ on $V$.
For $y\in V$, let $\tau_y$ be the \textbf{hitting time} of $y$:
\bq
\tau_y&\df& \inf\{t\geq 0\st X_t=y\}.\eq
\par
We are particularly interested in the functional
\bqn{F}
F\st \cL(G,\pi)\ni L&\mapsto & \sum_{x,y\in V}\pi(x)\pi(y)\EE^L_x[\tau_y],\eqn
where the subscript $x$ (respectively the superscript $L$) in the expectation indicates that $X$ is starting from $x$ (resp.\ is generated by $L$).
The probabilistic interpretation of $F(L)$ is the mean time to go from $x$ to $y$ for the Markov process generated by $L$, when $x$ and $y$ are sampled independently according to its invariant probability $\pi$.
\begin{rem}
The smaller the $F(L)$, the faster the underlying Markov process goes between the elements of $V$. It does not necessarily imply that the faster the time-marginal distributions go to equilibrium in large time (especially in the discrete time analogue). It is more related to the asymptotic behavior of the variance associated with the convergence of the empirical measures.
This point of view will not be investigated here.
\end{rem}\par
The quantity $F(L)$ also admits a nice spectral formulation:
for any $L\in\cL(G,\pi)$, let $\Lambda(L)$ be the spectrum of $-L$, removing the eigenvalue 0.
To take into account the possible multiplicities of the eigenvalues, $\Lambda(L)$ should be seen a \textbf{multiset}
(i.e.\ a eigenvalue of $-L$ of multiplicity $m$, appears $m$ times in $\Lambda(L)$).
By irreducibility of $L$, $\Lambda(L)$ is a priori a sub(multi)set of $\CC_+\df\{z\in \CC\st \Re(z)>0\}$ that is invariant under conjugation.
The \textbf{eigentime relation} asserts
\bqn{spectral1}
F(L)&=&\sum_{\lambda\in\Lambda(L)}\frac1\lambda.\eqn
The references  Cui and Mao \cite{MR2728712} and \cite{MR3442821} given in the discrete time setting also deal with the continuous time case.
The quantity $F(L)$ can also be written in terms of return times.
 Define for any $y\in V$,
\begin{eqnarray*}
\sigma &:=& \min\{t >0 : X_t \neq X_0\}, \\
T_y &:=& \min\{t \geq  \sigma : X_t = y\}.
\end{eqnarray*}
By irreducibility of $L$, we have the following eigentime identities (see Cui and Mao \cite{MR2728712}), for any $y\in V$:
\begin{eqnarray}
\sum_{x\in V}\pi(x)E_x[\tau_y]\ =\  \sum_{\lambda\in\Lambda(L)}\frac{1}{\lambda}\  =\ \sum_{x,z\in V}\pi(x)\pi(z)E_x[\tau_z], \label{eigen3} \\
\sum_x\pi(x)E_x[T_y] \ =\ 1 + \sum_{\lambda\in\Lambda(L)}\frac{1}{\lambda}\ =\ \sum_{x,z\in V}\pi(x)\pi(z)E_x[T_z], \label{eigen4}
\end{eqnarray}

Similarly to Subsection 2.3, our goal is to find the minimizers of $F$ on $\cL(G,\pi)$, or at least to deduce some information about them, since they correspond to the fastest normalized Markov processes compatible with $G$ with invariant distribution
 $\pi$. There is no loss of generality in imposing that $L$ is irreducible, because the functional $F$ is infinite for non-irreducible Markov generators admitting $\pi$ as invariant measure.\\

Consider next a continuous time $V$-valued \textit{controlled} Markov chain, denoted $(X_t)_{t \geq 0}$ again by abuse of notation, controlled by a control process $(Z_t)_{t \geq 0}$. The latter takes values in $\U_{X_t}$, where $\U_i := [0, \infty)^{|\V_i|}$, identified with the instantaneous transition rate of $X_t$. That is, as $\delta$ goes to $0_+$,
\begin{eqnarray*}
P(X_{t+\delta} = j | X_s, Z_s, s \leq t, X_t = i) &=& P(X_{t+\delta} = j | X_t = i, Z_t) \\
&=&\lt\{
\begin{array}{ll}
 Z_t(X_t,j)\delta + o(\delta)&\hbox{, if }j \in \V_{X_t}, \\
 -\sum_{j \in \V_{X_t}}Z_t(X_t,j)\delta + o(\delta)&\hbox{, if } j = X_t,\\
 0&\hbox{, if }  j \notin \V_{X_t}\cup\{X_t\},
\end{array}\rt.
\end{eqnarray*}
where we write $Z_t = (Z_t(X_t,j_1), \cdots, Z_t(X_t, j_{m_{X_t}}))$ for a suitable enumeration $(j_1, \cdots, j_{m_{X_t}})$ of $\V_{X_t}$.

For the remaining part of this subsection, we consider the case where $\pi=\upsilon$, the uniform measure on $V$.
The renormalization condition \eqref{normalisation} can be written in the form
\begin{equation}
\sum_{i \neq j}Z_t(i,j) \ =\ \vert V\vert. \label{normal}
\end{equation}

If for any $t\geq 0$, $Z_t$ is a function of $X_t$ alone, say $Z_t = r(X_t, \cdot) \in \U_{X_t}$, then $(X_t)_{t\geq 0}$ is a time-homogeneous Markov process with rate matrix $R = (r(i,j))_{i,j \in V}$, where we set $r(i,j) = 0$ for $j \notin \V_i$. Consider the problem of minimizing (\ref{eigen3}). As before, we augment the state process to the $V^*$-valued process $(\hat{X}_t)_{t\geq 0} = (X_t, A_t)_{t\geq 0}$, with the understanding that $A_t$ can change only when $X_t$ does and a transition of $X_t$ from $i$ to $j$ leads to a transition of $A_t$ to $A_t\backslash\{j\}$.  Consider the control problem of minimizing the  cost
\begin{equation}
E\left[\int_0^{\zeta}|A_t|dt \ \Big| \ X_0 = i_0, A_0 = V\backslash\{i_0\}\right], \label{cost'}
\end{equation}
for
$$\zeta := \{t \geq 0 : A_t = \emptyset\},$$
which is equivalent to (\ref{eigen3}), subject to the normalization constraint (\ref{normal}). The constraint (\ref{normal}) couples decisions across different states, so dynamic programming arguments cannot be directly applied. Therefore we modify the formulation for the time being, this modification will be dropped later. The modification is as follows.  Let $(a_i)_{i\in V}$ be scalars in $(0, \vert V\vert)$ such that $\sum_ia_i = \vert V\vert$. For state $i$, we restrict the rates to be from the set
$$\tilde{U}_i := \{r(i, j) : r(i,j) = 0 \ \forall \ j \notin \V_i, \  \sum_{j \in \V_i}r(i,j) = a_i\}.$$
Consider the value function
$$V(i, A) := \inf E\left[\int_0^{\zeta}|A_t|dt \ \Big| \ X_0 = i, A_0 = A\right],$$
where the infimum is over all admissible controls. The dynamic programming equation then is
\begin{equation}
\min_{r(i, \cdot) \in \tilde{U}_i}\left(|A| + \sum_{j \in \V_i, A}r(i,j)(V(j, A\backslash\{j\}) - V(i, A))\right) = 0, \ V(\cdot , \emptyset) \equiv 0. \label{DPc}
\end{equation}
Once again it is clear that the quantity being minimized is affine in the variables it is being minimized over and hence the optimum is attained for a deterministic choice of $r(i, \cdot)$ in the sense that $r(i,j)$ can be non-zero for at most one $j \in \V_i$. Thus the optimal path traces the nodes of $\G$ in a deterministic manner, visiting each of them at least once. This is true for any choice of $\{a_i\}$ and therefore true in general for the constraint $(\ref{normal})$. Unlike the discrete time case, this does not, however, mean that the trajectory is deterministic, because the sojourn time in each node is still random. It is clear that the cost for any such trajectory will be
$$E\left[\sum_{i \in \V}a_i^{-1}{N_i}\right],$$
where $N_i $ is the number of times the trajectory passed through $i$. For a given choice of $(a_i)_{i\in V}$, this is clearly minimized if $N_i = 1$ for all $i\in V$, which can be achieved  by tracing a Hamiltonian cycle if one exists. Optimizing next over the choice of $(a_i)_{i\in V}$ subject to (\ref{normal}), namely $\sum_{i\in V} a_i=\vert V\vert$, a simple induction argument shows that the choice $a_i = 1$, for all $i\in V$, is optimal.\\

As in the discrete case, a similar treatment is possible for the cost (\ref{eigen4}) or its equivalent
\begin{equation}
E\left[\int_0^{\zeta}|A_t|dt \ \Big| \ X_0 = i_0, A_0 = V\right], \label{cost'2}
\end{equation}
with the constraint (\ref{normal}).\\

We have proved:\\

\begin{theo}\label{theo1}  Minimum of either the cost (\ref{cost'}) or the cost (\ref{cost'2}) over all $V$-valued random processes compatible with $\G$ is attained by tracing a Hamiltonian cycle when one exists.
\end{theo}

 \par\me

 \section{Perturbation of Markov generators}

 \subsection{Introduction}

 We can expect the optimality of Hamiltonian cycle to persist under small perturbations of the Markov chains considered above. For specific classes of perturbations, such results were established in Ejov et al \cite{Ejov}. Here we establish a vastly more general result, first for the continuous time framework (which turns out to be more natural in some sense for the kind of techniques we employ) and then for the discrete case.\\

When $A\df (a_0, a_1, ..., a_{N-1})\in\cH$ and a positive probability measure $\pi$ on $V$ are fixed, the set $\cL(G_A,\pi)$ is reduced to a singleton, its element will be denoted $L_A$.
It is indeed given by
\bq
\fo x,y\in V,\qquad L_A(x,y)&=&\lt\{\begin{array}{ll}
\frac{1}{N\pi(x)}&\hbox{, if  $x= a_k$ and $y=a_{k+1}$ for some $k\in\ZZ_n$}\\
-\frac{1}{N\pi(x)}&\hbox{, if  $x=y$}\\
0&\hbox{, otherwise}\end{array}\rt.\eq

Theorem \ref{theo1} may seem a little deceptive: the fastest normalized Markov processes $X$ leaving invariant $\upsilon$, the uniform probability measure on $V$, follow a prescribed cyclic ordering of the states of $V$,
their randomness comes only from their waiting times, distributed according to the exponential law of intensity 1.
Such a derandomization of the successive points visited by $X$ is also valid for probability measures $\pi$ close to $\upsilon$:
\begin{theo}\label{theo2}
Assume that the graph $G$ is Hamiltonian.
Then there exists a neighborhood $\cN$ of $\upsilon$ in the set $\cP_+(V)$ of positive probability measures on $V$ (endowed with the topology inherited from that of $(0,1]^V$)
such that for any $\pi\in\cN$,
the set of minimizers of $F$ on $\cL(G,\pi)$ is exactly $\{L_A\st A\in \cH(G)\}$.
\end{theo}
\par
Nevertheless,  this result cannot be extended to all positive probability measures $\pi$, at least for the graphs which are not a Hamiltonian cycle, a situation where $\cL(G,\pi)$ is not reduced to a singleton, in particular, this requires $N\geq 3$.
\begin{theo}\label{theo3}
Assume that $G$ is not a Hamiltonian cycle.
Then there exist positive probability measures $\pi$ on $V$ such that none of the elements of
$\{L_A\st A\in \cH(G)\}$ is a minimizer of $F$ on $\cL(G,\pi)$.
\end{theo}
\par
Thus for some $(G,\pi)$, the minimizers of $F$ on $\cL(G,\pi)$ are (spatially) \textbf{hesitating} Markov processes: at some vertex, the next visited point is not chosen deterministically.
For a given Hamiltonian graph $G$ which is not reduced to a Hamiltonian cycle, it would be interesting to describe the probability measures $\pi$ leading to a transition between non-hesitating and hesitating minimizers. This issue remains open at present.\\

\subsection{Differentiation on $\cL(\pi)$}

This section introduces some elements of differential calculus on $\cL(\pi)$, which will be helpful in the proof of Theorem \ref{theo2}.
Here we will be working mainly with the complete graph $K_V$.\par\me
We begin by presenting a more analytical expression for the functional $F$.
 For $y\in V$, consider the function
\bq
f_y\st V\ni x&\mapsto & \frac{\un_{\{y\}}(x)}{\pi(y)}-1.\eq
Note that $\pi[f_y]=0$, so for any
$L\in \cL(\pi)$, by irreducibility, there exists a unique function $\varphi^L_y$ on $V$ satisfying the Poisson equation
\bqn{varphiy}
\lt\{\begin{array}{rcl}
L[\varphi^L_y]&=&f_y,\\
\varphi^L_y(y)&=&0.\end{array}\rt.\eqn
The following relation with the functional $F$ is well-known:
\begin{lem}\label{Fanalytic}
For any $L\in\cL(\pi)$ and any $x,y\in V$, we have
\bq
\varphi^L_y(x)&=& \EE^L_x[\tau_y],\eq
so that
\bq
F(L)&=&\sum_{y\in V}\pi(y)\pi[\varphi_y^L].\eq
\end{lem}
To simplify notation, from now on, we will remove the $L$ in the exponent of $\EE_x^L$ and $\varphi^L_y$, when the underlying generator $L$ is clear from the context.
\proof
Let us recall a simple argument, which will be used again in the sequel.
Through the martingale problem characterization of $X$, we have that for any given function $\varphi$ on $V$,
the process $(M_t)_{t\geq 0}$ defined by
\bq
\fo t\geq 0,\qquad  M_t&\df&\varphi(X_t)-\varphi(X_0)-\int_0^t L[\varphi](X_s)\, ds\eq
is a martingale.
In particular, for any stopping time $\tau$, the process $(M_{\tau\wedge t})_{t\geq 0}$ is also a martingale.
Thus, starting from $x\in V$, we get,
\bq
\EE_{x}[M_{\tau_y\wedge t}]&=&0\\
&=& \EE_x\lt[\varphi(X_{\tau\wedge t})-\varphi(X_0)-\int_0^{\tau\wedge t} L[\varphi](X_s)\, ds\rt]\\
&=& \EE_x\lt[\varphi(X_{\tau\wedge t})\rt]-\varphi(x)-\EE_x\lt[\int_0^{\tau\wedge t} L[\varphi](X_s)\, ds\rt].\eq
Since $\tau$ is a.s.\ finite and $\varphi(X_{\tau\wedge t}), t \geq 0$, uniformly integrable, we obtain, by letting $t$ go to infinity
\bq
 \EE_x\lt[\varphi(X_{\tau})\rt]-\varphi(x)-\EE_x\lt[\int_0^{\tau} L[\varphi](X_s)\, ds\rt]&=&0.\eq
For any $y\in V$, consider $\varphi\df \varphi_y$ and $\tau\df\tau_y$.
From \eqref{varphiy} and from the fact that $f_y(z)=-1$ for any $z\in V\setminus\{y\}$, we deduce
\bq
\varphi_y(x)&=& \EE_x[\tau_y].\eq
The last identity of the lemma comes from
\bq
\sum_{x\in V}\pi(x)\EE_x[\tau_y]
&=&\sum_{x\in V}\pi(x)\varphi_y(x)\\
&=&\pi[\varphi_y].\eq
\wwtbp
\par
Since we are looking for minimizers of $F$ on $\cL(\pi)$, it is natural to differentiate this functional.
Let $\bar\cL(\pi)$ be the convex set of  normalized Markov generators $L$  admitting $\pi$ for invariant probability.
The difference with $\cL(\pi)$ is that the elements of $\bar \cL(\pi)$ are not required to be irreducible.
For $L\in \cL(\pi)$, $\wi L\in \bar\cL(\pi)$ and $\epsilon\in[0,1)$, let $L_\epsilon\df (1-\epsilon) L+\epsilon \wi L\in\cL(\pi)$.
Define
\bq
D_{\wi L}F(L)&\df& \lim_{\epsilon\ri 0_+} \frac{F(L_\epsilon)-F(L)}{\epsilon}.\eq
In the proof of the following result, it will be shown that
 this limit exists.
\begin{lem}\label{analytic}
With the above notation, we have
\bq
D_{\wi L}F(L)&=&\sum_{y\in V} \pi(y)(\pi[\varphi_y]-\pi[\psi_y])\\
&=& F(L)-\sum_{y\in V} \pi(y)\pi[\psi_y].\eq
where
 $\psi_y$ is the unique solution of another Poisson equation
\bqn{psiy}
\lt\{\begin{array}{rcl}
L[\psi_y]&=&\wi L[\varphi_y],\\
\psi_y(y)&=&0\end{array}.\rt.\eqn
\end{lem}
\proof
Let $\cF_\pi$ stand for the space of functions $f$ on $V$ whose mean with respect to $\pi$ vanishes.
By restriction to $\cF_\pi$,  $L\in\cL(\pi)$ can be seen as an invertible endomorphism of $\cF_\pi$,
denote by $L^{-1}_{\vert \cF_\pi}$ its inverse. Similarly, for $\epsilon\in[0,1)$, let $L^{-1}_{\epsilon, \vert \cF_\pi}$
be the inverse of $L_\epsilon$ on $\cF_\pi$.
The mapping $[0,1)\ni \epsilon \mapsto L_\epsilon$ being analytical, the same is true for
$[0,1)\ni \epsilon \mapsto L^{-1}_{\epsilon, \vert \cF_\pi}$.
Since we have
\bq
\fo \epsilon\in [0,1),\,\fo y\in V,\qquad \varphi_{y}^{L_\epsilon}&=& L^{-1}_{\epsilon, \vert \cF_\pi}[f_y]-L^{-1}_{\epsilon, \vert \cF_\pi}[f_y](y),\eq
we deduce that the mapping
\bq
[0,1)\ni\epsilon&\mapsto & \varphi_{y}^{L_\epsilon}\eq
is analytical.
The same is true for $[0,1)\ni\epsilon\mapsto F(L_\epsilon)$, due to the equality
\bq
\fo \epsilon\in [0,1),\qquad F(L_\epsilon)&=&\sum_{y\in V}\pi(y)\pi[\varphi_{y}^{L_\epsilon}].\eq
In particular its derivative $D_{\wi L}F(L)$ exists and is equal to
$\sum_{y\in V}\pi(y)\pi[\varphi'_y]$, where
 $\varphi'_y$ is the derivative of $\varphi_y^{L_\epsilon}$ at $\epsilon=0$.
Differentiating the relation $L_\epsilon [\varphi_y^{L_\epsilon}]=f_y$, we get
\bq
(\wi L-L)[\varphi_y]+L[\varphi'_y]&=&0.\eq
Furthermore, we have that $\varphi'_y(y)=\pa_\epsilon \varphi^{L_\epsilon}_y(y)\vert_{\epsilon =0}=0$,
so that $\varphi_y-\varphi_y'$ satisfies the equation \eqref{psiy} and must be equal to $\psi_y$.
The claim then follows from the equality $\varphi'_y=\varphi_y-\psi_y$, for all $y\in V$.
\wwtbp
\par
In the above proof, we have seen that $[0,1)\ni\epsilon\mapsto F(L_\epsilon)$ is analytic, so we can differentiate it a second time at $\epsilon=0$.
Denote $D^2_{\wi L}F(L)=\pa^2_\epsilon F(L_\epsilon)\vert_{\epsilon=0}$.
\begin{lem}\label{D2wiL}
For $L\in \cL(\pi)$, $\wi L\in \bar\cL(\pi)$, we have
\bq
D^2_{\wi L}F(L)&=&\sum_{y\in V} \pi(y)(2\pi[\varphi_y]-4\pi[\psi_y]+2\pi[\psi'_y])\\
&=&4D_{\wi L}F(L)-2F(L)+2\sum_{y\in V}\pi(y)\pi[\psi_y'],\eq
where
 $\psi_y'$ is the unique solution of
\bqn{psiy1}
\lt\{\begin{array}{rcl}
L[\psi_y']&=&\wi L[\psi_y],\\
\psi'_y(y)&=&0\end{array}.\rt.\eqn
\end{lem}
\proof
For any $y\in V$,
denote $\varphi''_y$  the second derivative of $\varphi_y^{L_\epsilon}$ at $\epsilon=0$.
By differentiating twice the relation $L_\epsilon [\varphi_y^{L_\epsilon}]=f_y$ at $\epsilon=0$, we get
\bq
(\pa_\epsilon^2L_\epsilon)[\varphi_y]+2(\pa_\epsilon L_\epsilon)[\varphi'_y]+L[\varphi''_y]&=&0.\eq
namely, since $\pa_\epsilon^2L_\epsilon=0$,
\bq
L[\varphi''_y]&=&2(L-\wi L)[\varphi'_y]\\
&=&2(L-\wi L)[\varphi_y-\psi_y]\\
&=&2L[\varphi_y-\psi_y]-2\wi L[\varphi_y-\psi_y]\\
&=&2L[\varphi_y-\psi_y]-2\wi L[\varphi_y]+2\wi L[\psi_y]\\
&=&2L[\varphi_y-\psi_y]-2L[\psi_y]+2\wi L[\psi_y]\\
&=&L[2\varphi_y-4\psi_y]+2\wi L[\psi_y].
\eq
It follows that $\varphi''_y/2-\varphi_y+2\psi_y$ satisfies the first condition of equation \eqref{psiy1}.
It also vanishes at $y$, since $\varphi''_y(y)=0=\varphi_y(y)=\psi_y(y)$.
Thus we get that $\varphi''_y=2\varphi_y-4\psi_y+2\psi'_y$.
The announced result is now a consequence of the equality
\bq
D^2_{\wi L}F(L)&=&\sum_{y\in V} \pi(y)\pi[\varphi''_y].\eq
\wwtbp
\par
It  will be convenient to use these differentiations with respect to particular generators $\wi L\in\bar \cL(\pi)$.
A cycle $A$ in $V$ is a finite sequence $(a_0, a_1, ..., a_{n-1})$ of distinct elements of $V$, with $n\in\NN\setminus\{1\}$ (up to the identification with
$(a_k, a_{k+1}, ..., a_{k+n-1})$, for all $k\in\ZZ_n$). As with Hamiltonian cycles
(corresponding to $n=N$),
we will make the convention that $a_n=a_0$,
as the indices should be seen as elements of $\ZZ_n$. The set of all cycles is denoted by $\cA$.
For any $A\in\cA$, there is a unique element $L\in\bar\cL(\pi)$ such that
\bq \fo x,y\in V, \qquad L(x,y)>0 &\Leftrightarrow& \ex l\in \ZZ_n\st x= a_l\hbox{ and }y=a_{l+1}.\eq
It is indeed the generator, denoted $L_A$ in the sequel, given by
\bqn{LA}
\fo x,y\in V,\qquad L_A(x,y)&=&\lt\{\begin{array}{ll}
\frac{1}{n\pi(x)}&\hbox{, if  $x= a_l$ and $y=a_{l+1}$ for some $l\in\ZZ_n$},\\
-\frac{1}{n\pi(x)}&\hbox{, if  $x=y$},\\
0&\hbox{, otherwise}.\end{array}\rt.\eqn
\begin{lem}\label{der}
Let $A=(a_l)_{l\in\ZZ_n}\in\cA$ be given and for $y\in V$, consider the function $\psi_y$ defined by \eqref{psiy} with $\wi L=L_A$.
Then we have
\bqn{M2}
\fo x\in V,\qquad
\psi_y(x)&=&\frac1n\sum_{l\in \ZZ_n}(\varphi_y(a_{l+1})-\varphi_y(a_{l}))(\varphi_{a_l}(x)-\varphi_{a_l}(y)).\eqn
Furthermore, we get that
\bq
\sum_{y\in V}\pi(y)\pi[\psi_y]&=&-\frac1n\sum_{l\in \ZZ_n}\sum_{y\in V} \pi(y)(\varphi_y(a_{l+1})-\varphi_y(a_{l}))\varphi_{a_l}(y)\\
&=&\frac1n\lt(\sum_{l\in \ZZ_n}\frac12\EE_{a_{l+1}}[\tau_{a_{l}}^2]-\EE_\pi[\tau_{a_l}]\EE_{a_{l+1}}[\tau_{a_{l}}]\rt),\eq
 where  $\EE_\pi$ stands for the expectation relative to the initial distribution $\pi$ for $X$.\end{lem}
 \proof
 For any function $\varphi$ on $V$, we have
 \bq
 L_A[\varphi]&=&\sum_{l\in\ZZ_n} \frac{\varphi(a_{l+1})-\varphi(a_l)}{n\pi(a_l)}\un_{\{a_l\}}.\eq
 Let $\psi$ be a function such that $L[\psi]=L_A[\varphi]$.
Using the martingale problem as in the proof of Lemma~\ref{analytic}, we get for any $x,y\in V$,
 \bq
 \psi(y)-\psi(x)&=&\EE_x\lt[\int_0^{\tau_y}L[\psi](X_s)\, ds\rt]\\
 &=&\EE_x\lt[\int_0^{\tau_y}L_A[\varphi](X_s)\, ds\rt]\\
 &=&\sum_{l\in\ZZ_n} \frac{\varphi(a_{l+1})-\varphi(a_l)}{n}\EE_x\lt[\int_0^{\tau_y}\frac{\un_{\{a_l\}}}{\pi(a_l)}(X_s)\, ds\rt]\\
 &=& \sum_{l\in\ZZ_n} \frac{\varphi(a_{l+1})-\varphi(a_l)}{n}\EE_x\lt[\int_0^{\tau_y} 1+f_{a_l}(X_s)\, ds\rt]\\
 &=& \sum_{l\in\ZZ_n} \frac{\varphi(a_{l+1})-\varphi(a_l)}{n}\lt( \varphi_y(x)+\EE_x\lt[\int_0^{\tau_y} f_{a_l}(X_s)\, ds\rt]\rt).\eq
 Taking into account that $L[\varphi_{a_l}]=f_{a_l}$,
 we deduce that
 \bq
 \EE_x\lt[\int_0^{\tau_y} f_{a_l}(X_s)\, ds\rt]&=&\varphi_{a_l}(y)-\varphi_{a_l}(x),\eq
 so that
 \bq
  \psi(y)-\psi(x)&=&\sum_{l\in\ZZ_n} \frac{\varphi(a_{l+1})-\varphi(a_l)}{n}\lt( \varphi_y(x)+\varphi_{a_l}(y)-\varphi_{a_l}(x)\rt).\eq
 Note that
 \bq
 \sum_{l\in\ZZ_n} \frac{\varphi(a_{l+1})-\varphi(a_l)}{n} \varphi_y(x)&=& \frac{\varphi_y(x)}{n}\sum_{l\in\ZZ_n} \varphi(a_{l+1})-\varphi(a_l)\\
 &=&0.\eq
 Thus
 \bq
  \psi(y)-\psi(x)&=&\frac1n\sum_{l\in\ZZ_n} (\varphi(a_{l+1})-\varphi(a_l))\lt( \varphi_{a_l}(y)-\varphi_{a_l}(x)\rt).\eq
 Considering for $y\in V$ the functions $\varphi=\varphi_{y}$ and $\psi=\psi_y$ and recalling that $\psi_y(y)=0$,
gives the first relation of the lemma.
Integrating this relation with respect to $\pi$ in $x$,
we get
\bq
\pi[\psi_y]&=&\frac1n\sum_{l\in \ZZ_n}(\varphi_y(a_{l+1})-\varphi_y(a_{l}))(\pi[\varphi_{a_l}]-\varphi_{a_l}(y)).\eq
A well-known result (recall \eqref{eigen1} or see e.g.\ the book of Aldous and Fill \cite{aldous-fill-2014}) asserts that
the quantity $\sum_{y\in V} \pi(y)\varphi_y(x)$ does not depend on $x\in V$.
It follows that
\bqn{Aldous}
\sum_{y\in V}\pi(y)(\varphi_y(a_{l+1})-\varphi_y(a_{l}))\pi[\varphi_{a_l}]&=&0\eqn
and hence
\bq
\sum_{y\in V}\pi(y)\pi[\psi_y]&=&-\frac1n\sum_{y\in V}\pi(y)\sum_{l\in \ZZ_n}(\varphi_y(a_{l+1})-\varphi_y(a_{l}))\varphi_{a_l}(y),\eq
which is the second equality of the lemma.
For any $l\in \ZZ_n$, let $\phi_{a_l}$ be the function defined by:
\bqn{phial}
\fo x\in V,\qquad \phi_{a_l}(x)&=&\sum_{y\in V} \pi(y)\varphi_{a_l}(y)(\varphi_y(x)-\varphi_y(a_l)).\eqn
We have $\sum_{y\in V}\pi(y)\pi[\psi_y]=-\frac1n\sum_{l\in \ZZ_n}\phi_{a_l}(a_{l+1})$.
To compute $\phi_{a_l}$, note that $\phi_{a_l}(a_l)=0$ and that
\bq
L[\phi_{a_l}]&=& \sum_{y\in V} \pi(y)\varphi_{a_l}(y)L[\varphi_y]\\
&=&\sum_{y\in V} \pi(y)\varphi_{a_l}(y)f_y\\
&=&\sum_{y\in V} \pi(y)\varphi_{a_l}(y)\lt(\frac{\un_{\{y\}}}{\pi(y)}-1\rt)\\
&=&\varphi_{a_l}-\pi[\varphi_{a_l}].\eq
This observation leads us to resort once again to the martingale problem, to get
for any $x\in V$,
\bq
\phi_{a_l}(a_l)&=&\phi_{a_l}(x)+\EE_x\lt[\int_0^{\tau_{a_l}}\varphi_{a_l}(X_s)-\pi[\varphi_{a_l}]\, ds\rt]\\
&=& \phi_{a_l}(x)+\EE_x\lt[\int_0^{\tau_{a_l}}\varphi_{a_l}(X_s)\, ds\rt]-\pi[\varphi_{a_l}]\EE_x[\tau_{a_l}]\\
&=&\phi_{a_l}(x)-\pi[\varphi_{a_l}]\varphi_{a_l}(x)+\EE_x\lt[\int_0^{\tau_{a_l}}\varphi_{a_l}(X_s)\, ds\rt]\\
&=&\phi_{a_l}(x)-\pi[\varphi_{a_l}]\varphi_{a_l}(x)+\frac12\EE_x[\tau_{a_l}^2]\eq
according to Lemma \ref{tau2} below.
Recalling that $\phi_{a_l}(a_l)=0$, we get
\bqn{phil}
\phi_{a_l}(x)&=&\pi[\varphi_{a_l}]\varphi_{a_l}(x)-\frac12\EE_x[\tau_{a_l}^2]\eqn
and this leads immediately to the last equality of the lemma.\wwtbp
In the previous proof, we needed the following result.
\begin{lem}\label{tau2}
For any $x,y\in V$,  we have
\bq
\EE_x\lt[\int_0^{\tau_{y}}\varphi_{y}(X_s)\, ds\rt]&=& \frac12\EE_x[\tau_{y}^2].\eq
 \end{lem}
 \proof
 Coming back to the probabilistic interpretation of $\varphi_y$, we get
 \bq
 \EE_x\lt[\int_0^{\tau_{y}}\varphi_{y}(X_s)\, ds\rt]&=&\int_0^{+\iy}\EE_x\lt[\un_{\{s\leq \tau_y\}}\EE_{X_s}[\tau_y]\rt]\, ds\\
 &=&\int_0^{+\iy}\int_0^{+\iy}\EE_x\lt[\un_{\{s\leq \tau_y\}}\EE_{X_s}[\un_{\{t\leq \tau_y\}}]\rt]\, ds\,dt\\
 &=&\int_0^{+\iy}\int_0^{+\iy}\EE_x\lt[\un_{\{s\leq \tau_y\}}\EE[\un_{\{t\leq \tau_y\circ \theta_s\}}\vert \sigma(X_u\st u\in[0,s])]\rt]\, ds\,dt\\
 &=&\int_0^{+\iy}\int_0^{+\iy}\EE_x\lt[\un_{\{s\leq \tau_y\}}\un_{\{t\leq \tau_y\circ \theta_s\}}\rt]\, ds\,dt\\
&=&\int_0^{+\iy}\int_0^{+\iy}\EE_x\lt[\un_{\{s+t\leq \tau_y\}}\rt]\, ds\,dt,
 \eq
 where we used the Markov property and
 where $\theta_s$ is the shift by time $s$  of the trajectories of $X$.
Using the Fubini theorem, we get
\bq
\int_0^{+\iy}\int_0^{+\iy}\EE_x\lt[\un_{\{s+t\leq \tau_y\}}\rt]\, ds\,dt&=&
\EE_x\lt[\int_0^{+\iy}\int_0^{+\iy}\un_{\{s+t\leq \tau_y\}}\, ds\,dt\rt]\\
&=&\frac12\EE_x\lt[\int_0^{+\iy}\int_0^{+\iy}\un_{\{s\leq \tau_y\}}\un_{\{t\leq \tau_y\}}\, ds\,dt\rt]\\
&=&\frac12\EE_x\lt[\lt(\int_0^{\tau_y}\, ds\rt)^2\rt]\\
&=&\frac12\EE_x[\tau_{y}^2]\eq
\wwtbp
\par
Before treating  the second derivative in a similar way, let us
present two remarks about the quantities entering Lemma \ref{der}. We believe they will be
relevant for further study of
the
minimizers of the mapping $F$ on $\cL(\pi)$.
\par
Define the following quantities, associated with a given $L\in\cL(\pi)$:
\bqn{HAL}
\nonumber\fo x,y\in V,\qquad h_L(x,y)&\df& \frac12\EE_{y}[\tau_{x}^2]-\EE_\pi[\tau_{x}]\EE_{y}[\tau_{x}]\\
\fo A=(a_1, ..., a_n)\in\cA,\qquad H_A(L)&=&\frac1n\sum_{l\in \ZZ_n} h_L(a_l,a_{l+1})
\eqn
Lemma \ref{der} can be rewritten under the form
\bqn{DFH}
\fo A=(a_1, ..., a_n)\in\cA,\qquad
D_AF(L)&=&F(L)-H_A(L)
\eqn
where $D_AF(L)$ is a short hand for $D_{L_A}F(L)$.
\par
Let us say that a cycle $A=(a_1, ..., a_n)\in\cA$ is below the generator $L$, if
\bq
\fo l\in\ZZ_n,\qquad L(a_l,a_{l+1})>0\eq
and denote by $\cA(L)$ the set of cycles below $L$.
Then we have:
\begin{lem}
Assume that $L\in \cL(\pi)$ is a minimizer of $F$ on $\cL(\pi)$. Then,
\bq
\fo A\in\cA(L),\qquad H_A(L)&=&F(L),\\
\fo A\in \cA\setminus\cA(L),\qquad H_A(L)&\leq &F(L).\eq
In particular, we get
\bq F(L)&=&\max_{A\in\cA} H_A(L).\eq
\end{lem}
\proof
Consider  a minimizer $L\in \cL(\pi)$ of $F$ on $\cL(\pi)$ and $A\in\cA(L)$.
Then for $\epsilon\in\RR$ small enough, $(1-\epsilon)L+\epsilon L_A$ remains a Markov generator
and belongs to $\cL(\pi)$. Differentiating $F(L_\epsilon)$ at $\epsilon=0$, we thus get that
$D_AF(L)=0$, which implies $H_A(L)=F(L)$.
For $A\in \cA\setminus\cA(L)$, the operator $(1-\epsilon)L+\epsilon L_A$ is not Markovian for $\epsilon <0$.
So $D_AF(L)$ only corresponds  to the right derivative of $F(L_\epsilon)$ at $\epsilon=0_+$.
The minimizing assumption on $L$  implies that $D_AF(L)\geq 0$, namely
$H_A(L)\leq F(L)$.
The last identity of the lemma is an immediate consequence of the previous observations and of the fact that there exists
at least one cycle below $L$, by irreducibility.\wwtbp
\par
Next we mention a spectral relation satisfied by the quantities $(h_L(x,y))_{x,y\in V}$, reminiscent of \eqref{spectral1}.
Indeed, it is proved in a similar way, as  will become clear from the following proof where the arguments for \eqref{spectral1} will be recalled.
\par
\begin{lem}
For any $L\in\cL(\pi)$, he have
\bqn{spectral2}
\sum_{x,y\in V}\pi(x)\pi(y)h_L(x,y)&=&\sum_{\lambda\in\Lambda(L)}\frac1{\lambda^2}.\eqn
\end{lem}
\proof
As in the proof of Lemma \ref{Fanalytic}, let
$\cF_\pi$ stand for the space of functions $f$ on $V$ whose mean with respect to $\pi$ vanishes
and denote by $\Pi$ the orthogonal projection from $\LL^2(\pi)$ to $\cF_\pi$:
\bq
\fo f\in \LL^2(\pi),\qquad \Pi[f]&=&f-\pi[f].\eq
Let $(g_y)_{y\in V}$ be an orthonormal basis of $\LL^2(\pi)$ and $R$ be any endomorphism of $\cF_\pi$.
We have seen in Lemma 6 of \cite{MR3442821} that
\bq
\trace(R)&=&\sum_{y\in V} \pi[\Pi[g_y]R[\Pi[g_y]]].\eq
In \cite{MR3442821}, we considered the orthonormal basis given by
\bq
\fo y\in V,\qquad g_y&\df& \frac{\un_{\{y\}}}{\sqrt{\pi(y)}}
\eq
and the operator $L^{-1}_{\vert \cF_\pi}$ defined in Lemma \ref{Fanalytic}, in order to conclude \eqref{spectral1},
taking into account the fact that $\Pi[g_y]=\sqrt{\pi(y)}f_y$, for all $y\in V$, and that
$\trace(L^{-1}_{\vert \cF_\pi})=\sum_{\lambda\in\Lambda(L)}\frac1\lambda$.
\par
To prove \eqref{spectral2}, we use $R=(L^{-1}_{\vert \cF_\pi})^2$.
Remark that for any $y\in V$,
\bq
(L^{-1}_{\vert \cF_\pi})^2[f_y]
&=& L^{-1}_{\vert \cF_\pi}[\varphi_y-\pi[\varphi_y]]\\
&=& \phi_y-\pi[\phi_y]\eq
where $\phi_y$ is the unique solution of
\bqn{phi}
\lt\{\begin{array}{rcl}
L[\phi_y]&=&\varphi_y-\pi[\varphi_y]\\
\phi_y(y)&=&0.\end{array}\rt.\eqn
(This notation agrees with that introduced in \eqref{phial}).
Thus we get
\bq
\sum_{\lambda\in\Lambda(L)}\frac1{\lambda^2}&=&
\trace( (L^{-1}_{\vert \cF_\pi})^2)\\
&=&\sum_{y\in V} \pi[\Pi[g_y](L^{-1}_{\vert \cF_\pi})^2[\Pi[g_y]]]\\
&=&\sum_{y\in V} \pi(y) \pi[f_y(L^{-1}_{\vert \cF_\pi})^2[f_y]]\\
&=&\sum_{y\in V} \pi(y) \pi[f_y(\phi_y-\pi[\phi_y])]\\
&=&\sum_{y\in V} \pi(y) \pi[f_y\phi_y]\\
&=&\sum_{y\in V} \pi(y)\phi_y(y)-\pi(y) \pi[\phi_y]\\
&=&-\sum_{y\in V} \pi(y) \pi[\phi_y]\\
&=&-\sum_{x,y\in V}\pi[x]\pi[y] \phi_y(x).
\eq
In the proof of Lemma \ref{der} (see \eqref{phil}), it was shown that
\bq
\fo x,y\in V,\qquad \phi_x(y)&=&-h_L(x,y),\eq
which leads immediately to \eqref{spectral2}.
\wwtbp
\par
Lemma \ref{der} can be extended to the second derivative presented in Lemma \ref{D2wiL}, by computing similarly the function $\psi'_y$ defined by \eqref{psiy1} with $\wi L=L_A$, for fixed $A=(a_l)_{l\in\ZZ_n}\in\cA$ and $y\in V$.
For our purposes, it is convenient to consider a generalization of this situation.
Given another cycle $A'=(a'_l)_{l\in\ZZ_{n'}}\in\cA$, consider the equation in the function $\Psi_y$:
\bqn{psiy2}
\lt\{\begin{array}{rcl}
L[\Psi_y]&=&L_{A'}[\psi_y],\\
\Psi_y(y)&=&0,\end{array}\rt.\eqn
where $\psi_y$ is still associated to $L$, $A$ and $y$ as in Lemma \ref{der}.
Of course, when $A'=A$,  we recover $\Psi_y=\psi'_y$.
\begin{lem}\label{derdeux}
For $A=(a_l)_{l\in\ZZ_n}\in\cA$,  $A'=(a'_l)_{l\in\ZZ_{n'}}\in\cA$ and $y\in V$  given as above, consider the function $\Psi_y$ defined by \eqref{psiy2}.
Then we have, for any $ x\in V$,
\bqn{M3}
\Psi_y(x)&=&\frac1{nn'}\sum_{l\in \ZZ_n, k\in\ZZ_{n'}}(\varphi_y(a_{l+1})-\varphi_y(a_{l}))(\varphi_{a_l}(a'_{k+1})-\varphi_{a_l}(a'_k))(\varphi_{a'_k}(x)-\varphi_{a'_k}(y)).\eqn
Furthermore, we get that
\bqn{der2A}
\nonumber\lefteqn{\sum_{y\in V}\pi(y)\pi[\Psi_y]}\\&=&\frac1{nn'}\sum_{l\in \ZZ_n, k\in\ZZ_{n'}}(h_L(a'_k,a_{l+1})-h_L(a'_k,a_{l}))(\varphi_{a_l}(a'_{k+1})-\varphi_{a_l}(a'_k))\\
\nonumber&=&\frac1{nn'}\sum_{l\in \ZZ_n, k\in\ZZ_{n'}}\lt(\frac12(\EE_{a_{l+1}}[\tau_{a'_k}^2]-\EE_{a_{l}}[\tau_{a'_k}^2])-\EE_\pi[\tau_{a'_k}](\EE_{a_{l+1}}[\tau_{a'_k}]-\EE_{a_{l}}[\tau_{a'_k}]).
\rt)\lt( \EE_{a'_{k+1}}[\tau_{a_l}]-  \EE_{a'_{k}}[\tau_{a_l}] \rt)
\eqn
\end{lem}
 \proof
From Lemma \ref{der}, we have
\bq
L_{A'}[\psi_y]&=&\frac1n\sum_{l\in \ZZ_n}(\varphi_y(a_{l+1})-\varphi_y(a_{l}))L_{A'}[\varphi_{a_l}]\\
&=&
\frac1n\sum_{l\in \ZZ_n}(\varphi_y(a_{l+1})-\varphi_y(a_{l}))\sum_{k\in\ZZ_{n'}}
 \frac{\varphi_{a_l}(a'_{k+1})-\varphi_{a_l}(a'_k)}{n'\pi(a'_k)}\un_{\{a'_k\}}\\
 &=&
\frac1{nn'}\sum_{l\in \ZZ_n}(\varphi_y(a_{l+1})-\varphi_y(a_{l}))\sum_{k\in\ZZ_{n'}}
 (\varphi_{a_l}(a'_{k+1})-\varphi_{a_l}(a'_k))\frac{\un_{\{a'_k\}}}{\pi(a'_k)}\\
 &=&\frac1{nn'}\sum_{l\in \ZZ_n}(\varphi_y(a_{l+1})-\varphi_y(a_{l}))\sum_{k\in\ZZ_{n'}}
 (\varphi_{a_l}(a'_{k+1})-\varphi_{a_l}(a'_k))f_{a'_k},\eq
 where we used that for any $l\in\ZZ_n$,
 \bq
 \sum_{k\in\ZZ_{n'}}
 {\varphi_{a_l}(a'_{k+1})-\varphi_{a_l}(a'_k)}&=&0.\eq
 Thus, denoting
 \bq
 \xi_y&\df&\frac1{nn'}\sum_{l\in \ZZ_n}(\varphi_y(a_{l+1})-\varphi_y(a_{l}))\sum_{k\in\ZZ_{n'}}
 (\varphi_{a_l}(a'_{k+1})-\varphi_{a_l}(a'_k))(\varphi_{a'_k}-\varphi_{a'_k}(y)),\eq
 we get that $L[\xi_y]=L_{A'}[\psi_y]$ and $\xi_y(y)=0$.
 It follows that $\Psi_y=\xi_y$, as announced.
\par
We deduce that
\bq
\pi[\psi'_y]&=&\frac1{nn'}\sum_{l\in \ZZ_n}(\varphi_y(a_{l+1})-\varphi_y(a_{l}))\sum_{k\in\ZZ_{n'}}
 (\varphi_{a_l}(a'_{k+1})-\varphi_{a_l}(a'_k))(\pi[\varphi_{a'_k}]-\varphi_{a'_k}(y))\eq
 and
 \bqn{tech}
\nonumber \sum_{y\in V}\pi(y)\pi[\psi'_y]&=&\frac1{nn'}\sum_{y\in V}\pi(y)\sum_{l\in \ZZ_n, k\in\ZZ_{n'}}(\varphi_y(a_{l+1})-\varphi_y(a_{l}))
 (\varphi_{a_l}(a'_{k+1})-\varphi_{a_l}(a'_k))(\pi[\varphi_{a'_k}]-\varphi_{a'_k}(y))\\
 &=&-\frac1{nn'}\sum_{y\in V}\sum_{l\in \ZZ_n, k\in\ZZ_{n'}}\pi(y)(\varphi_y(a_{l+1})-\varphi_y(a_{l}))
 (\varphi_{a_l}(a'_{k+1})-\varphi_{a_l}(a'_k))\varphi_{a'_k}(y),\eqn
 where we used again (recall \eqref{Aldous}) that
 \bq
\sum_{y\in V} \pi(y)(\varphi_y(a_{l+1})-\varphi_y(a_{l}))&=&0.
 \eq
 Remember also (cf.\ \eqref{phial}) that
 \bq
 \fo x\in V,\qquad \sum_{y\in V}\pi(y)(\varphi_y(x)-\varphi_y(a'_{k}))
\varphi_{a'_k}(y)&=&\phi_{a'_k}(x)\\
&=&-h_L(a'_k,x).\eq
Thus substituting in \eqref{tech}
 \bq
 \varphi_y(a_{l+1})-\varphi_y(a_{l})&=&\varphi_y(a_{l+1})-\varphi_y(a'_{k})-(\varphi_y(a_{l})-\varphi_y(a'_{k}))\eq
 we deduce \eqref{der2A}.
 The last equality of the lemma is obtained by expressing $h_L$ and $\varphi_x$, for $x\in V$, in terms of expectation of hitting times.
 \wwtbp
\par
Denote by $H_{A',A}(L)$ the expression given by \eqref{der2A}.
Considering the case $A'=A$, Lemma \ref{D2wiL} leads to
\bq
D_{A,A}F(L)&=&
2F(L)+4H_A(L)+2H_{A,A}(L)\eq
where $D_{A,A}F(L)$ is a shorthand for $D^2_{L_A}F(L)$.
But the importance of Lemma \ref{derdeux}, is because,
if we define for any $A,A'\in \cA$,
$D_{A',A}F(L)\df D_{A'}(D_AF(L))$,
then we get
\bq
D_{A',A}F(L)&=&
2(F(L)+H_A(L)+H_{A'}(L)+H_{A',A}(L)).\eq\par
The previous expressions for the differentiations up to order 2 with respect to Markov generators associated to cycles can be extended to general Markov
generators from $\bar \cL(\pi)$.
To go in this direction, we need to recall a simple result:
\begin{lem}\label{extremal}
The extremal points of the convex set $\bar\cL(\pi)$
are exactly the generators $L_A$ for $A\in\cA$.
\end{lem}
As a consequence, any $L\in\bar \cL(\pi)$ can be decomposed into a barycentric sum
\bq
 L&=&\sum_{A\in \cA} p(A)L_A,\eq
 where $p$ is a probability measure on $\cA$.
For an extensive discussion of such decompositions, see the book of Kalpazidou \cite{MR2226353}.
Note that the above decomposition is not unique in general, because $\bar\cL(\pi)$ is not a simplex for $N\geq 3$.
For instance, the generator
\bq
L&\df& \frac12\lt(\begin{array}{ccc}
-2&1&1\\
1&-2&1\\
1&1&-2\end{array}\rt)\eq
of the simple random walk on $\ZZ_3$ can be written in the form
$L=\frac12L_{(0,1,2)}+\frac12L_{(0,2,1)}$ and $L=\frac13L_{(0,1)}+\frac13L_{(1,2)}+\frac13L_{(2,0)}$.
\par
Nevertheless, given $\wi L, \wit L\in \bar \cL(\pi)$, decompose them into
\bq
 \wi L&=&\sum_{A\in \cA} \wi p(A)L_A,\\
  \wit L&=&\sum_{A\in \cA} \wit p(A)L_A,\eq
 where $\wi p, \wit p$ are probability measures on $\cA$.
 Then we get for any $L\in \cL(\pi)$.
 \bq
 D_{\wi L}F(L)&=&\sum_{A\in \cA} \wi p(A)D_AF(L),\\
 D_{\wit L}D_{\wi L}F(L)&=&\sum_{A,A'\in \cA} \wi p(A)\wit p(A') D_{A,A'}F(L).\eq
 It follows that we can  write
 \bq
  D_{\wi L}F(L)&=&F(L)-H_{\wi L}(L),\\
 D_{\wit L}D_{\wi L}F(L)&=&2(F(L)+H_{\wi L}(L)+H_{\wit L}(L)+H_{\wi L,\wit L}(L)),\eq
 where
 \bq
 H_{\wi L}(L)&=&\sum_{x\not =y} \pi(x)L(x,y) h_L(x,y)\\
 H_{\wit L, \wi L}(L)&=&\sum_{x\not= y, \, x'\not =y'} \pi(x')\wit L(x',y')\pi(x)\wi L(x,y) (h(x',y)-h(x',x))(\varphi_x(y')-\varphi_x(x'))\eq
 (definitions which conform to \eqref{HAL} and \eqref{der2A} when $\wi L=L_A$ and $\wit L=L_{A'}$).
 \par\sm
 In view of \eqref{M2} and \eqref{M3},
 the following quantity seems to play an important role in bounding the derivatives:
 \bq
 M(L)&\df&\max_{y,x,x'\in V}\lve \varphi_y(x)-\varphi_y(x')\rve\\
 &=&\max_{y,x\in V} \varphi_y(x).\eq
 \begin{pro}\label{borneder}
 We have for any $L\in\cL(\pi)$ and $\wi L,\,\wit L\in \bar\cL(\pi)$,
 \bq
 F(L)&\leq & M(L),\\
 \lve D_{\wi L}F(L)\rve&\leq & M(L)+M(L)^2,\\
 \lve D_{\wit L}D_{\wi L}F(L)\rve &\leq & 2(M(L)+M(L)^2+M(L)^3).\eq
 \end{pro}
 \proof
 The first bound is obvious. For the second, note that
 \eqref{M2} can be extended to the solution of \eqref{psiy} for general $\wi L\in\bar \cL(\pi)$: we get
\bq
\fo y,x\in V,\qquad
\psi_y(x)&=&\sum_{z\not=z'\in V}\pi(z)\wi L(z,z')(\varphi_y(z')-\varphi_y(z))(\varphi_{z}(x)-\varphi_{z}(y)).\eq
Taking into account the renormalization of $\wi L$, it follows that for any $y\in V$, we have for the supremum norm:
\bq
\lVe \psi_y\rVe_\iy&\leq & M(L)^2\eq\par
For the third bound of the lemma, note that \eqref{M3} can also be extended to $\Psi_y$ for given $y\in V$, which is the solution of
\bq
\lt\{\begin{array}{rcl}
L[\Psi_y]&=&\wit L [\psi_y]\\
\Psi_y(y)&=&0,\end{array}\rt.\eq
where $\psi_y$ is the solution of \eqref{psiy}. It follows that for any $x,y\in V$,
\bq
\Psi_y(x)&=&\sum_{u\not= v, \, u'\not =v'} \pi(u')\wit L(u',v')\pi(u)\wi L(u,v)(\varphi_y(v)-\varphi_y(u))(\varphi_{u}(v')-\varphi_{u}(u'))(\varphi_{u'}(x)-\varphi_{u'}(y)).\eq
Thus
\bq
\lVe \Psi_y\rVe_\iy&\leq & M(L)^3.\eq
\wwtbp\par
A natural question is how to upper bound $M(L)$. A first answer is to use the operator norm $\vvv \cdot \vvv_{\iy\ri\iy}$ from $\LL^\iy(\pi)$ to $\LL^\iy(\pi)$ with the operator $L^{-1}_{\vert \cF_\pi}$
introduced in Lemma \ref{Fanalytic}:
\bq
M(L)&=&\max_{y\in V}\lVe \varphi_y\rVe_\iy\\
&\leq & \max_{y\in V}\lVe \varphi_y-\pi[\varphi_y]\rVe_\iy+\max_{y\in V}\pi[\varphi_y]\\
&\leq & \vvv L^{-1}_{\vert \cF_\pi}\vvv_{\iy\ri\iy}\max_{y\in V}\lVe f_y\rVe_\iy+\max_{y\in V}(1-\pi(y))M(L)\\
&\leq & \vvv L^{-1}_{\vert \cF_\pi}\vvv_{\iy\ri\iy}\frac{1}{\pi_\wedge}+(1-\pi_\wedge)M(L)\eq
where $\pi_\wedge\df \min_{x\in V}\pi(x)$.
It follows that
\bq
M(L)&\leq & \frac{\vvv L^{-1}_{\vert \cF_\pi}\vvv_{\iy\ri\iy}}{\pi^2_\wedge}.\eq

But the norm $\vvv L^{-1}_{\vert \cF_\pi}\vvv_{\iy\ri\iy}$ does not seem so easy to evaluate. One can instead resort to  the operator norm  from $\LL^2(\pi)$ to $\LL^2(\pi)$ as follows. Denoting $I$ the identity operator on $\cF_\pi$, we have as above
\bq
M(L)&\leq & \vvv I\vvv_{2\ri \iy}\vvv L^{-1}_{\vert \cF_\pi}\vvv_{2\ri 2}\max_{y\in V}\lVe f_y\rVe_2+(1-\pi_\wedge)M(L)\\
&\leq &\frac1{\sqrt{\pi_\wedge}}\vvv L^{-1}_{\vert \cF_\pi}\vvv_{2\ri 2}\max_{y\in V}\sqrt{\frac1{\pi(y)}+1-\pi(y)}+(1-\pi_\wedge)M(L)\\
&\leq &\frac1{\sqrt{\pi_\wedge}}\vvv L^{-1}_{\vert \cF_\pi}\vvv_{2\ri 2}\sqrt{1+\frac{1}{\pi_\wedge}}+(1-\pi_\wedge)M(L)\\
&\leq & \vvv L^{-1}_{\vert \cF_\pi}\vvv_{2\ri 2}\frac{\sqrt{2}}{\pi_\wedge}+(1-\pi_\wedge)M(L).\eq
As a consequence, we get
\bq
M(L)&\leq & \frac{\sqrt{2}\vvv L^{-1}_{\vert \cF_\pi}\vvv_{2\ri 2}}{\pi^2_\wedge}.\eq
This expression is advantageous when $L$ is reversible with respect to $\pi$,
since in this situation, $\vvv L^{-1}_{\vert \cF_\pi}\vvv_{2\ri 2}=1/\lambda$, where $\lambda$ is the spectral gap of $L$, namely the smallest element of $\Lambda(L)$ (which is then  in $(0,+\iy)$).
Nevertheless, since we are interested in $F(L)$, note there is a simple comparison:
\bqn{MLF}
M(L)&\leq & \frac{F(L)}{\pi_\wedge^2}.\eqn
\par\me
We now concentrate on the case $\pi=\upsilon$, the uniform measure and $L=L_A$, with $A$ a Hamiltonian cycle.
The following result will be crucial in the proof of Theorem \ref{theo2}.
\begin{pro}\label{grander}
For any $A\in\cH$ and $\wi A\in \cA\setminus\{A\}$, we have on $\cL(\upsilon)$,
\bq
D_{\wi A} F(L_A)&\geq & \frac{N-1}{2N}.\eq
\end{pro}
\proof
There is no loss of generality in assuming that $V=\ZZ_N$ and that $A=(0,1,2, ..., N-1)$.
To simplify the notation, let us write $L=L_A$.
By invariance of $L$ and $\upsilon$ through the rotations $\ZZ_N\ni x\mapsto x+y\in\ZZ_N$ for any fixed $y\in\ZZ_N$,
it follows that the quantity $\EE_\upsilon[\tau_x]$ does not depend on the choice of $x\in\ZZ_N$.
It is then necessarily equal to $F(L)$.
Furthermore, since under $L$, the Markov process waits an exponential time before adding 1 to the current state,
we get that for any $x,y\in\ZZ_N$, $\EE_x[\tau_y]=\rho(x,y)$, where
\bq
\fo x,y\in\ZZ_N,\qquad \rho(x,y)&\df&\min\{n\in\ZZ_+\st y=x+n\}.\eq
It follows easily that $F(L)=(N-1)/2$ (for an alternative proof, see Corollary \ref{cor1} in the next section).
Thus we get that
\bq
\fo x,y\in\ZZ_N,\qquad
h_L(x,y)&=& \frac12\lt(\EE_y[\tau_x^2]-(N-1)\EE_y[\tau_x]\rt)\\
&=&\frac12\lt(\EE_y[\tau_x^2]-(N-1)\rho(y,x)\rt).\eq
Since under $\PP_y$, $\tau_x$ is a sum of $\rho(y,x)$ independent exponential random variables of parameter 1,
we compute that
\bq
\EE_y[\tau_x^2]&=&\EE_y[\tau_x]^2+\rho(y,x)\eq
(because for any exponential random variable $\cE$, we have $\EE[\cE^2]=2\EE[\cE]^2$).
Thus we get that for any $x,y\in\ZZ_N$, $h_L(x,y)=h_N(\rho(y,x))$, where
\bq
h_N\st [0,N-1]\ni r&\mapsto &\frac12\lt(r^2-(N-2)r\rt).\eq
This function $h_N$ is decreasing on $[0,(N-2)/2]$, increasing on $[(N-2)/2, N-1]$ and we have
$h_N(0)=0<(N-1)/2=h_N(N-1)$ .\par
Thus from the definition \eqref{HAL},
we get that
\bq
\fo \wi A\in\cA,\qquad H_{\wi A}(L)&\leq & h_N(N-1)\\
&=&H_{A}(L).\eq
More precisely, with $\wi A\df(a_0, a_1, ...,  a_{ n})\in\cA$, we get, except if for any $l\in\ZZ_n$, $h_N(\rho(a_l,a_{l+1}))=h_N(N-1)$,
\bq
H_{\wi A}(L)&\leq& \frac{n-1}{n}h_N(N-1)+\frac1n\max\{h_N(0),h_N(N-2)\}\\
&=&\frac{n-1}{n}H_A(L)\\
&\leq &\frac{N-1}{N}H_A(L),
\eq
where in the equality, we used that $h_N(0)=h_N(N-2)=0$ and that $h_N(N-1)=H_A(L)$, according to \eqref{HAL}.
But if for any $l\in\ZZ_n$, we have $h_N(\rho(a_l,a_{l+1}))=h_N(N-1)$,
it means that $a_{l+1}=a_l+1$. Since this must be true for all $l\in\ZZ_n$, it follows that $n=N$
and that $\wi A$ must be of the form $(k,{k+1}, ..., {k+N-1})$, for some $k\in\ZZ_N$, namely, it is Hamiltonian.\par
From \eqref{DFH}, we obtain,
\bq
\fo \wi A\in\cA\setminus\{A\}, \qquad D_{\wi A} F(L)
&=&F(L)-H_{\wi A}(L)\\ &\geq &F(L)-\frac{N-1}{N}H_{\wi A}(L)\\
&= &D_{A} F(L)+\frac1NH_A(L)\\
&=&\frac{N-1}{2N}\eq
due to $D_{A} F(L)=0$, because $L$ is not modified by modifying it in the direction of the cycle $A$.
\wwtbp\par\me
Above we worked with the complete graph $K_V$ and the associated set of Markov generators $\cL(\pi)$. But all the previous considerations can be extended to the
case of $\cL(G,\pi)$, where the graph $G$ is as in the introduction. The only difference is that $\cA$ has to be replaced by $\cA(G)$, the set of cycles using only edges from $E$.
For instance, Lemma \ref{extremal} has to be replaced by
\begin{lem}
The extremal points of the convex set $\bar\cL(G,\pi)$ (the set of  normalized Markov generators $L$, compatible with $G$  and admitting $\pi$ for invariant probability)
are exactly the generators $L_A$ for $A\in\cA(G)$.
\end{lem}

\subsection{Perturbations of the uniform probability measure}

Our main goal here is to show Theorems \ref{theo2} and \ref{theo3}. Their proofs are respectively based on small and large perturbations  of the uniform probability measure $\upsilon$.\par\me
First we check that all Hamiltonian cycles have the same speed in $\cL(\pi)$, as was announced in the introduction in the discrete time setting and for the uniform distribution $\upsilon$, but this is true more generally.
\begin{lem}\label{Hamiltonian}
Let $A=(a_0, ..., a_{N-1})\in\cH$ be a Hamiltonian cycle,
 we have
\bq
F(L_A)&=&\frac{N}2\sum_{x\not=y}\pi(x)\pi(y).\eq
In particular this quantity does not depend on the choice of the Hamiltonian cycle $A$.
\end{lem}
\proof
The generator $L_A$ can be represented by the matrix
\bq
\lt(
\begin{array}{cccccc}
-\frac{1}{N\pi(a_0)}&\frac{1}{N\pi(a_0)}&0&\cdots&\cdots&0\\
0&-\frac{1}{N\pi(a_1)}&\frac{1}{N\pi(a_1)}&0&\cdots&0\\
 & & & & & \\
 & &\cdots&\cdots& & \\
  & & & & & \\
  0&\cdots&\cdots& 0&-\frac{1}{N\pi(a_{N-2})}&\frac{1}{N\pi(a_{N-2})}\\
\frac{1}{N\pi(a_{N-1})}&0&\cdots&\cdots& 0&-\frac{1}{N\pi(a_{N-1})}
\end{array}\rt).
\eq
It follows that the polynomial in $X$ given by
\bq
P(X)&\df&\det \lt(
\begin{array}{cccccc}
X-\frac{1}{N\pi(a_0)}&\frac{1}{N\pi(a_0)}&0&\cdots&\cdots&0\\
0&X-\frac{1}{N\pi(a_1)}&\frac{1}{N\pi(a_1)}&0&\cdots&0\\
 & & & & & \\
 & &\cdots&\cdots& & \\
  & & & & & \\
  0&\cdots&\cdots&0& X-\frac{1}{N\pi(a_{N-2})}&\frac{1}{N\pi(a_{N-2})}\\
\frac{1}{N\pi(a_{N-1})}&0&\cdots&\cdots& 0&X-\frac{1}{N\pi(a_{N-1})}
\end{array}\rt),
\eq
is equal to $X\prod_{\lambda\in\Lambda(L_A)}(X-\lambda)$.
Expanding the latter expression into $X(\alpha_0+\alpha_1X+\cdots +\alpha_{N-1}X^{N-1})$, we get that
\bq
\sum_{\lambda\in\Lambda(L_A)}\frac1{\lambda}&=&-\frac{\alpha_1}{\alpha_0}\eq
This is indeed a consequence of
\bq
\alpha_0&=&(-1)^{N-1}\prod_{m\in\lin N-1\rin } \theta_m\\
\alpha_1&=&(-1)^{N-2}\sum_{k\in\lin N-1\rin}\prod_{m\in\lin N-1\rin \setminus\{k\}} \theta_m
\eq
where   $\Lambda(L_A)$ is parametrized as the multiset consisting of the $\theta_m$, for $m\in\lin N-1\rin\df\{1,2, ..., N-1\}$.
\par
On another hand, we compute directly from the definition of $P(X)$, by expanding the determinant, that
\bq
P(X)&=&\prod_{l\in\ZZ_N}\lt(X-\frac1{N\pi(a_l)}\rt)-\prod_{l\in\ZZ_N}\lt(-\frac1{N\pi(a_l)}\rt)\eq
It follows that
\bq
\alpha_0&=&\sum_{k\in\ZZ_N}\prod_{m\in\ZZ_N\setminus\{k\}}\lt(-\frac1{N\pi(a_m)}\rt)\\
\alpha_1&=&\frac12\sum_{k\not= l\in\ZZ_N}\prod_{m\in\ZZ_N\setminus\{k,l\}}\lt(-\frac1{N\pi(a_m)}\rt)\eq
(the factor $1/2$ is due to the fact that the couple $(k,l)$ also appears as $(l,k)$).
Multiplying the numerator and the denominator by
$\prod_{m\in\ZZ_N}\lt(-{N\pi(a_m)}\rt)$, we get that
\bq
-\frac{\alpha_1}{\alpha_0}&=&\frac{\frac{N^2}2\sum_{k\not= l\in\ZZ_N}\pi(a_k)\pi(a_l)}{N\sum_{m\in \ZZ_N}\pi(a_m)}\eq
and this leads to the announced result.\wwtbp
In particular, for $\pi=\upsilon$, the uniform probability measure on $V$, we get:
\begin{cor}\label{cor1}
For $\pi=\upsilon$, we have for any Hamiltonian cycle $A$,
\bq
F(L_A)&=&\frac{N-1}2.\eq
\end{cor}
\par
The next result is the crucial step in the proof of Theorem \ref{theo2}.
For its statement, introduce for any $A\in\cH$ and $\epsilon \in(0,1)$,
\bqn{NAe}
\cN_{A,\epsilon}&\df&
\{L=(1-t)L_A+t\wi L\st t\in[0,\epsilon)\hbox{ and }\wi L\in\bar\cL(\pi)\}
\eqn
This set is a neighborhood of $L_A$ in $\cL(\pi)$ and observe that we would have ended
with the same set if we had required in this definition that $\wi L$ belong to the convex hull generated by
the $L_{\wi A}$, for $\wi A\in \cA\setminus\{A\}$.\par
Define
\bq
\epsilon_1(N,\pi_\wedge)&\df& \pi_\wedge^4\ln\lt(1+\frac{1}{N\pi_\wedge^2}\rt),\\
\epsilon_2(\pi_\wedge)&\df&\frac{1}{56} \pi_\wedge^{12}, \\
\epsilon(N,\pi_\wedge)&\df& \epsilon_1(N,\pi_\wedge)\wedge \epsilon_2(\pi_\wedge).
\eq
\begin{lem}\label{crucial}
For $N\geq 2$ and any $A\in\cH$, $L_A$ is the unique minimizer of $F$ over $\cN_{A,\epsilon(N,\pi_\wedge)}$.
\end{lem}
\proof
Assume that for some given $A\in\cH$, $L_A$ is not the unique minimizer of $F$ over $\cN_{A,\epsilon(N,\pi_\wedge)}$.
Then we can find $t\in(0,\epsilon(N,\pi_\wedge))$ and a probability $p$ on $\cA\setminus\{A\}$, such that
$F(L_t)\leq F(L_A)$, with
\bq
L_t&\df& (1-t)L_A+t\wi L,\\
\wi L&\df& \sum_{\wi A\in \cA\setminus\{A\}}p(\wi A)L_{\wi A}.
\eq
Applying Taylor-Lagrange formula to the function
$[0,t]\ni s\mapsto F(L_s)$,
we get there exists $s\in[0,t]$ such that
\bq
F(L_t)&=&F(L_A)+tD_{\wi L}F(L_A)+\frac{t^2}{2} D^2_{\wi L}F(L_s).\eq
Taking into account Propositions \ref{borneder} and \ref{grander} and \eqref{MLF}, we obtain
\bqn{FLtFLs}
F(L_t)&\geq &F(L_A)+t\frac{N-1}{2N}-t^2\lt(\frac{F(L_s)}{\pi_\wedge^2}+\frac{F(L_s)^2}{\pi_\wedge^4}+\frac{F(L_s)^3}{\pi_\wedge^6}\rt).\eqn
To evaluate $F(L_s)$,
note that for $s\in(0,t)$,
\bq
\pa_s F(L_s)&=& D_{\wi L}F(L_s)\\
&\leq & M(L_s)+M(L_s)^2\\
&\leq & \frac{F(L_s)}{\pi_\wedge^2}+\frac{F(L_s)^2}{\pi_\wedge^4}.\eq
Classical computations show that if a $\cC^1$ function $f\st [0,t]\ri (0,+\iy)$ satisfies
$\pa_s f(s)\leq af(s)+bf^2(s)$ for all $s\in[0,t]$, where $a,b>0$, then assuming
$f(0)\exp(bt)< f(0)+b/a$,
we get
\bq
\fo s\in [0,t], \qquad f(s)&\leq & \frac{bf(0)\exp(bt)}{b+af(0)(1-\exp(bt))}.\eq
In particular, if
\bqn{condfs}
\exp(bt)&<& 1+b/(2af(0)),
\eqn
 then
 \bq
 \fo s\in [0,t], \qquad f(s)&\leq &2\lt(f(0)+\frac{b}{2a}\rt).\eq
Let us apply this observation with the mapping $[0,t]\ni s\mapsto F(L_s)$ and $a\df1/\pi_\wedge^2$, $b=1/\pi_\wedge^4$.
Since \bq
F(L_0)&=&F(L_A)\\
&=&\frac{N}2\sum_{x\not=y\in V}\pi(x)\pi(y)\\
&\leq &\frac{N}2, \eq
we get that condition \eqref{condfs} is satisfied, due to the definition of $\epsilon_1(N,\pi_\wedge)$ and
to the fact that $t\in(0,\epsilon_1(N,\pi_\wedge))$.
 It follows that,
\bq
\fo s\in[0,t],\qquad F(L_s)
&\leq & N+\frac1{\pi_\wedge^2}\\
&\leq & \frac2{\pi_\wedge^2}.
\eq
since $\pi_\wedge\leq 1/N$.
Substituting this bound in \eqref{FLtFLs}, we deduce that
\bq
F(L_t)&\geq &F(L_A)+t\frac{N-1}{2N}-t^2\lt(\frac{2}{\pi_\wedge^4}+\frac{4}{\pi_\wedge^8}+\frac{8}{\pi_\wedge^{12}}\rt)\\
&\geq &F(L_A)+t\frac{1}{4}-\frac{14}{\pi_\wedge^{12}}t^2
\eq
The r.h.s.\ is strictly larger than $F(L_A)$
if $t<\epsilon_2(\pi_\wedge)$ and this is in contradiction with our initial assumption.
\wwtbp
\par
Denote for any $\pi\in\cP_{+}(V)$,
\bqn{Fwedge}
F_\wedge(\pi)&\df&\inf\{F(L)\st L\in \cL(\pi)\}.\eqn
Another ingredient in the proof of Theorem \ref{theo2} is:
\begin{lem}\label{contFw}
The mapping $\cP_+(V)\ni\pi\mapsto F_\wedge(\pi)$ is continuous.
\end{lem}
\proof
Let $\cL$ be the set of irreducible and normalized Markov generators (so that $\cL=\sqcup_{\pi\in\cP_+(V)}\cL(\pi)$), endowed with the topology inherited from
$\RR^{V^2}$.
The functional $F$ is defined on $\cL$ and \eqref{spectral1} is valid on $\cL$.
As a consequence, $F$ is continuous on $\cL$.
Indeed, if $(L_n)_{n\in\NN}$ is a sequence of elements of $\cL$ converging to $L\in\cL$, then
according to Paragraph 5 of Chapter 2 of Kato \cite{MR1335452}, we have $\lim_{n\ri\iy} \Lambda(L_n)=\Lambda(L)$ and so
$\lim_{n\ri\iy} F(L_n)=F(L)$.
Next consider a sequence $(\pi_n)_{n\in\NN}$ of elements from $\cP_+(V)$ converging to $\pi\in\cP_+(V)$ and such that the sequence $(F_\wedge(\pi_n))_{n\in\NN}$ admits a limit.
For all $n\in\NN$, let $L_n$ be an element from $\cL(\pi_n)$ such that
\bq
F_\wedge(\pi_n)\ \leq \ F(L_n)\ \leq \ F_\wedge(\pi_n)+\frac1n\eq
Due to the normalization condition and to the belonging of $\pi$ to $\cP_+(V)$,
we can extract a subsequence (still denoted $(L_n)_{n\in\NN}$ below) from $(L_n)_{n\in\NN}$ converging to some generator $L$.
It is clear that $L$ is normalized and that $\pi$ is invariant for $L$.
Let us check that $L$ is irreducible. Fix $x\in V$.
For any $n\in\NN$, let $X^{(n)}\df (X^{(n)}_t)_{t\geq 0}$ be a Markov process starting from $x$ and whose generator is $L_n$.
It is not difficult to deduce from the corresponding martingale problems, that $X^{(n)}$ converges in law (with respect to the Skorokhod topology) to a Markov process starting from $x$
and whose generator is $L$. Thus for any $y\in V$ and $T\geq 0$,
\bq
\lim_{n\ri\iy}\EE_x[T\wedge\tau^{(n)}_y]&=&\EE_x[T\wedge\tau_y]\eq
(with an obvious notation). It follows that
\bq
\EE_x[T\wedge\tau_y]&\leq &\liminf_{n\ri\iy} \EE_x[\tau^{(n)}_y]\\&\leq &
\liminf_{n\ri\iy}\frac1{\pi_n(x)\pi_n(y)}F(L_n)\\
&=&\frac1{\pi(x)\pi(y)}\liminf_{n\ri\iy}F_\wedge(\pi_n)\\
&\leq & \frac{N}{2\pi(x)\pi(y)}\eq
according to Lemma \ref{Hamiltonian}.
Letting $T$ go to infinity, we get
$\EE_x[\tau_y]\leq N/(2\pi_\wedge^2)$.
This bound, valid for all $x,y\in V$, implies that $L$ is irreducible and thus $L\in\cL(\pi)$.
Furthermore, the above arguments show that
\bq
\lim_{n\ri\iy} F_\wedge(\pi_n)&=&\lim_{n\ri\iy} F(L_n)\\
&=&F(L)\\&\geq & F_\wedge(\pi).\eq
So $F_\wedge$ is lower continuous on $\cP_+(V)$. By considering the sequence $(\pi_n)_{n\in\NN}$ identically equal to $\pi$,
we also get that the infimum defining $F_\wedge(\pi)$ is attained.
\par
To show that $F_\wedge$ is upper continuous on $\cP_+(V)$,
let again $(\pi_n)_{n\in\NN}$ be a sequence of elements from $\cP_+(V)$ converging to some $\pi\in\cP_+(V)$  and such that the sequence $(F_\wedge(\pi_n))_{n\in\NN}$ admits a limit.
According to the previous remark, there exists $L\in\cL(\pi)$ such that
$F(L)=F_\wedge(\pi)$.
For any $n\in\NN$, consider the matrix $\wi L_n$ given by
\bq
\fo x,y\in V,\qquad \wi L_n(x,y)&\df&\frac{\pi(x)}{\pi_n(x)}L(x,y).\eq
It is immediate to prove that $\wi L_n$ is an irreducible Markov generator leaving $\pi_n$ invariant.
But it may not be normalized, so let $\kappa_n>0$ be such that $L_n\df\kappa_n \wi L_n$ belongs to $\cL(\pi_n)$.
There is no difficulty in checking that $L_n$ converges to $L$  and thus that
$\lim_{n\ri\iy} F(L_n)=F(L)=F_\wedge(\pi)$. Thus passing into the limit in
$F(L_n)\geq F_\wedge(\pi_n)$, we deduce that
\bq
F_\wedge(\pi)&\geq &\lim_{n\ri\iy} F_\wedge(\pi_n)\eq
as desired.
\wwtbp
\par
With all these ingredients, we can now come to the
\prooff{Proof of Theorem \ref{theo2}}
Note that it is sufficient to consider the case where $G$ is the complete graph over $V$, since $F_\wedge(\pi)\leq \min\{F(L)\st L\in\cL(G,\pi)\}$, for any graph $G$
and positive probability measure $\pi$ on $V$.\par
The main argument is by contradiction. Assuming that the statement of Theorem \ref{theo2} is not true,
we can find a sequence $(\pi_n)_{n\in\NN}$ converging to $\upsilon$, such that for all $n\in\NN$,
there exists $L_n\in\cL(\pi_n)\setminus\{L_{\pi_n,A}\st A\in \cH\}$ with $F(L_n)=F_\wedge(\pi_n)$.
(Here we have included  $\pi_n$ in the index of $L_{\pi_n,A}$ to underscore the fact that this generator, associated to a Hamiltonian cycle $A$,
also depends on the underlying invariant probability $\pi_n$.)
As seen in the proof of Lemma \ref{contFw}, a subsequence (still denoted $(L_n)_{n\in\NN}$)  converging toward some $L\in \cL(\upsilon)$ can be extracted from $(L_n)_{n\in\NN}$.
We furthermore have
\bq
\lim_{n\ri\iy} F(L_n)&=&F(L)\eq
and by Lemma \ref{contFw}
\bq
\lim_{n\ri\iy} F_\wedge(\pi_n)&=&F_\wedge(\upsilon).\eq
It follows that $F(L)=F_\wedge(\upsilon)$.
From Theorem \ref{theo1}, we deduce that there exists $A\in\cH$ such that $L=L_{\upsilon, A}$.
Using again the fact that
\bqn{piup}\lim_{n\ri\iy} \pi_n&=&\upsilon,\eqn
 we get that
$\lim_{n\ri\iy} L_{\pi_n,A}=L_{\upsilon, A}$ and thus
\bqn{LL0}\lim_{n\ri\iy}(L_n- L_{\pi_n,A})&=&0.\eqn
Consider $r\df\min_{n\in\NN} \pi_{n,\wedge}$, which is positive due to \eqref{piup},
and let $\epsilon\df \epsilon(N,r)$, with the notation introduced before Lemma \ref{crucial}.
From \eqref{LL0}, we deduce that for $n\in\NN$ large enough,
$L_n$ belongs to $\cN(\pi_n,A,\epsilon)$, defined as in \eqref{NAe}, with $\pi$ replaced by $\pi_n$.
Then Lemma \ref{crucial} asserts that $L_n=L_{\pi_n,A}$, because $L_n$ is a minimizer of $L$ over $\cL(\pi_n)$.
This is in contradiction with our initial assumption.
\wwtbp
\par
To finish this section,  we consider large  perturbations of the uniform probability measure $\upsilon$.
\prooff{Proof of Theorem \ref{theo3}}
Let $G=(V,E)$ be a finite oriented connected graph which is not a Hamiltonian cycle.
Then we can find a cycle $A\df(a_0, a_1, ..., a_{n-1})\in\cA(G)$ with $n<\card(V)$.
Denote $\wi V\df\{a_0, a_1, ..., a_{n-1}\}$ and $\wit V\df V\setminus\wi V$.
By the strong connectivity of $G$, we can find a subset $\wit E$ of oriented edges from $E$,
such that $\card(\wit E)=\card(\wit V)$ and for any $x\in \wit V$ we can find exactly one $y\in V$
with $(x,y)\in \wit E$.
Putting together the edges from $A$ and those from $\wit E$, we get a graph $\breve G$ on $V$ looking like the following picture, where the cycle is oriented clockwise and the trees are oriented toward the cycle.
\begin{figure}[H]
 \centering
  \includegraphics[width=6cm]{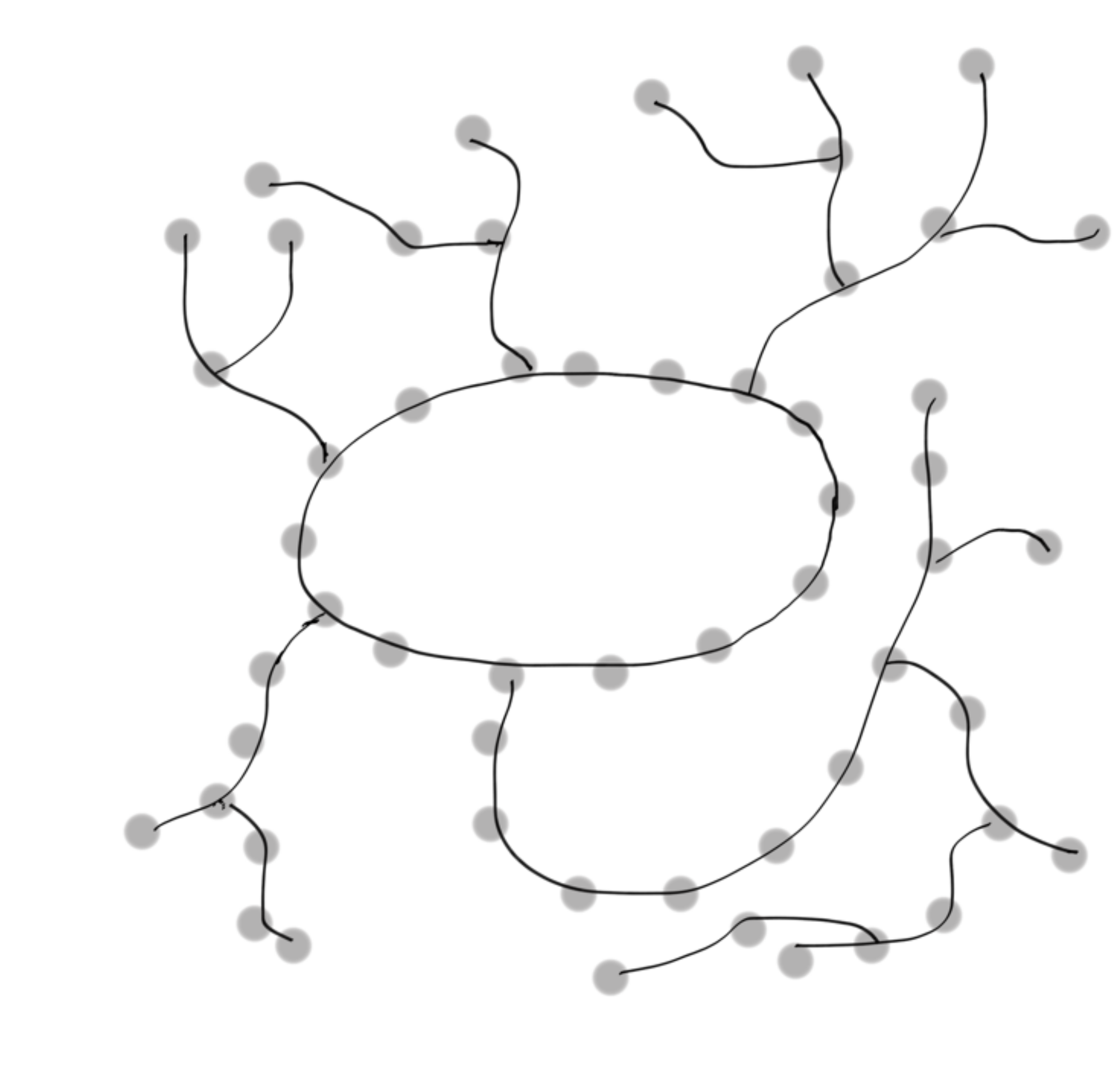}
  \caption{the graph $\breve G$}
 \end{figure}
For $r>0$, consider the Markov generator $L_r$ defined by
\bq
\fo x\not=y\in V,\qquad
L_r(x,y)&\df&\lt\{\begin{array}{ll}
1&\hbox{, if there exists $l\in\lin0,n-1\rin$ such that $x=a_l$ and $y=a_{l+1}$}\\
r&\hbox{, if $(x,y)\in\wit E$}\\
0&\hbox{, otherwise}.
\end{array}\rt.
\eq
This generator is not irreducible, since it does not allow the chain to go from the cycle $A$ to $\wit V$.
Nevertheless, its unique invariant probability measure $\pi$ is the uniform probability measure on $\wi V$.
The generator $L_r$ then satisfies an \textbf{extended normalization condition}, in the sense that
\bq
\sum_{x\not=y\in V} \pi(x)L_r(x,y)&=&1\eq
The interest in $L_r$ is because it is easy to find its eigenvalues:
\bq
\Lambda(L_r)&=&\Lambda(\wi L_A)\sqcup\{ r[\vert \wit V\vert]\}\eq
where $\wi L_A$ is the generator corresponding to the Hamiltonian cycle given by $A$ on $\wi V$ and $\{r[\vert \wit V\vert]\}$ is the multiset
consisting of the value $r$ with the multiplicity $\vert \wit V\vert$.
This identity is an  immediate consequence of following decomposition of $L_r$, where all the elements of $\wi V$ have been put before those of $\wit V$ and where the elements of $\wit V$ have been ordered so that the (oriented) distance to $\wi V$ is non-decreasing (in particular the last element corresponds to a leaf of $\breve G$):
\bq
L_r&=&
\lt(
\begin{array}{cc}
\wi L_A &0\\
C&D
\end{array}\rt).
\eq
In the r.h.s., the $\wit V\times\wit V$ matrix $D$  is sub-diagonal and its diagonal consists only of $-r$.
Formula~\eqref{spectral1} enables to extend the functional $F$ to $L_r$ and we get
\bq
F(L_r)&=&F(\wi L_A)+\frac{\vert V\vert}{r}.\eq
In particular, it follows that
\bq
\lim_{r\ri+\iy} F(L_r)\ =\ F(\wi L_A)
\ =\ \frac{n-1}{2}\ <\ F(L_H)\eq
for any Hamiltonian cycle $H\in\cH(G)$,
where we used twice Corollary \ref{cor1}.
From now on, we fix $r>0$ large enough, so that
\bqn{FAFH}
F(L_r)&<& F(L_H)\eqn
for any Hamiltonian cycle $H\in\cH(G)$.
\par
For any $\epsilon>0$, consider the Markov generator
\bq
 L_{r,\epsilon}&\df& Z_{r,\epsilon}^{-1}(L_r+\epsilon L_G)\eq
 where \\
 $\bullet$ the Markov generator $L_G$ is defined by
 \bq
 \fo x\not =y\in V,\qquad
 L_G(x,y)&\df& \lt\{\begin{array}{ll}
1&\hbox{, if $(x,y)\in E$}\\
0&\hbox{, otherwise}
\end{array}\rt.
\eq
$\bullet$ the constant $Z_{r,\epsilon}>0$ is such that $L_{r,\epsilon}$ is normalized (this is possible because $L_r+\epsilon L_G$
is irreducible on $V$).\par
For $r,\epsilon>0$, denote $\pi_{r,\epsilon}$ the invariant probability measure of $L_{r,\epsilon}$.
It is clear that as $\epsilon$ goes to $0_+$, $\pi_{r,\epsilon}$ converges toward $\pi$.
It follows that
\bq
\lim_{\epsilon\ri 0_+}Z_{r,\epsilon}&=&1\\
\lim_{\epsilon\ri 0_+}L_{r,\epsilon}&=&L_r\eq
From the general theory of perturbation of spectra of finite operators (see e.g.\ the beginning of the second chapter of the book of Kato \cite{MR1335452}), we have
\bq
\lim_{\epsilon\ri 0_+} F(L_{r,\epsilon})&=&F(L_r).\eq
Taking into account \eqref{FAFH}, we can thus find $\epsilon>0$ small enough so that
\bq
F(L_{r,\epsilon})&<& F(L_H)\eq
for any Hamiltonian cycle $H\in\cH(G)$.
Namely the probability measure $\pi_{r,\epsilon}$ satisfies the statement of Theorem \ref{theo3}.
One would have remarked that this  probability measure $\pi_{r,\epsilon}$ is quite far away from $\upsilon$, because it gives
very small weight to the elements of $\wit V$.
\wwtbp

\section{The discrete time framework}

Here we  discuss the links between the search of the fastest continuous-time Markov processes
with the analogous problem in discrete time.\par\me
Let a graph $G=(V,E)$ and a positive probability measure $\pi$ on $V$ be fixed and denote
by $ \cK(G,\pi)$  the set of irreducible Markov kernels $K$ on $V$ whose permitted transitions are edges from $E$ (plus self-loops, i.e., the possibility to stay at the same place) and leaving $\pi$ invariant, namely satisfying $\pi K=\pi$.
For any $K\in\cK(G,\pi)$, let $\fX\df(\fX_n)_{n\in\ZZ_+}$ be a Markov chain whose transitions are dictated by $K$.
For any $y\in V$, recall (see \eqref{hitting}) that
\bq
\uptau_y&\df&\inf\{n\in \ZZ_+\st \fX_n=y\}.\eq
On $\cK(G,\pi)$, we consider the functional $\mathfrak{F}$ defined by
\bq
\fo K\in  \cK(G,\pi), \qquad\mathfrak{F}(K)&\df&\sum_{x,y\in V} \pi(x)\pi(y)\EE_x[\uptau_y],\eq
where subscript $x$ in the expectation indicates that $\fX$ is starting from $x\in V$.
\par
To any $K\in\cK(G,\pi)$, we associate $\Theta(K)$ the multiset consisting of the spectrum of $K$, removing the eigenvalue 1 (of multiplicity 1). It is a priori a sub(multi)set of the closed unit disk centered at 0 of $\CC$ left invariant by conjugation.
Analogously to the continuous-time situation, we have the eigentime relation
\bq
\fo K\in\cK(G,\pi),\qquad \fF(K)&=&\sum_{\theta\in \Theta(K)}\frac1{1-\theta}.\eq
\par
To any $L\in\cL(G,\pi)$, associate
\bq
l&\df&\max\{L(x)\st x\in V\},\\
K&\df&I+\frac{L}{l}.
\eq
It is immediate to check that $K\in\cK(G,\pi)$.
Furthermore, we have $\Theta(K)=1-\Lambda(L)/l$, so that
\bqn{fFF}
\fF(K)&=&l F(L).\eqn
Taking into account that
\bq
l\ \geq \ \sum_{x\in V} \pi(x) L(x)\ =\ 1,\eq
it follows that $\fF(K)\geq F(L)$. We will denote $\Phi\st \cK(G,\pi)\ri \cL(G,\pi)$ the mapping  $L \to K$ defined above.
\par
Conversely, to any $K\in \cK(G,\pi)$, associate
\bq
k&\df& \frac{1}{\sum_{x\in V} \pi(x)(1-K(x,x))},\\
L&=&k(K-I).\eq
It is immediate to check that $L\in \cL(G,\pi)$.
Furthermore, we get $\Lambda(L)=k(1-\Theta(K))$
and it follows that
\bq
F(L)&=&\fF(K)/k.\eq
Taking into account that
\bq
k\ \geq \  \frac{1}{\sum_{x\in V} \pi(x)}\ =1,\eq
we get that $F(L)\leq \fF(K)$. Denote $\Psi\st \cL(G,\pi)\ri \cK(G,\pi)$ the mapping  $L \to K$ as above.\par
\begin{rem}
The mappings $\Phi$ and $\Psi$ are not inverse of each other, because the image of $\cL(G,\pi)$ by $\Phi$ is included
into $\cK_0(G,\pi)\df\{ K\in \cK(G,\pi)\st \ex x\in V\hbox{ with } K(x,x)=0\}$.
Nevertheless, we have that $\Phi$ and $\Psi_0$ are inverse of each other, where $\Psi_0$ is the restriction of $\Psi$ to $\cK_0(G,\pi)$.
\par When one is looking for the minimal value of $\fF$ on $\cK(G,\pi)$, one can restrict  attention to $\cK_0(G,\pi)$, because
\bq
\min\{\fF(K)\st K\in\cK(G,\pi)\}&=&\min\{\fF(K)\st K\in\cK_0(G,\pi)\}.\eq
Indeed, for any $K\in\cK(G,\pi)$, there exist a unique $\wi K\in\cK_0(G,\pi)$ and $\alpha\in [0,1)$ such that $K=(1-\alpha)\wi K+\alpha I$.
Then we get $\Theta(K)=(1-\alpha)\Theta(\wi K)+\alpha$, i.e.\ $\Theta-1=(1-\alpha)(\Theta(\wi K)-1)$.
This implies that
\bq
\fF(\wi K)&=&(1-\alpha)\fF(K)\\
&\leq & \fF(K).\eq
\end{rem}
\par
As in \eqref{Fwedge}, denote
\bq
F_\wedge(G,\pi)&\df&\inf\{F(L)\st L\in \cL(G,\pi)\},\\
\fF_\wedge(G,\pi)&\df&\inf\{\fF(K)\st K\in \cK(G,\pi)\}.\eq
From the above considerations, we deduce:
\begin{pro}\label{contdiscr}
We always have
\bq
F_\wedge(G,\pi)&\leq &\fF_\wedge(G,\pi)\eq
(in particular, when looking for the fastest Markov chain in the sense we have defined, it is preferable to resort to continuous time rather than to discrete time).
\par
Furthermore, assume that there is a minimizer $L\in \cL(G,\pi)$ of $F$ such that $L(x)$ does not depend on $x\in V$ (it is then equal to 1).
Then $F_\wedge(G,\pi)=\fF_\wedge(G,\pi)$.
\end{pro}
\proof
Consider $K\in\cK(G,\pi)$. We have seen that
\bq
\fF(K)&\geq & F(\Psi(K))\\&\geq &F_\wedge(G,\pi),\eq
so taking the infimum over
$K\in\cK(G,\pi)$, we get the first bound.
\par
Conversely, if $L\in \cL(G,\pi)$ is a minimizer of $F$ whose diagonal is constant, then $l=1$ in \eqref{fFF}, namely
$\fF(\Phi(L))=F(L)=F_\wedge(G,\pi)$.
From the previous inequality, it follows that $\Phi(L)$ is indeed a minimizer of $\fF$ on $\cK(G,\pi)$
and we conclude that $F_\wedge(G,\pi)=\fF_\wedge(G,\pi)$.
\wwtbp
\par
In association with Theorem \ref{theo1}, the above proposition also enables us to recover the result of
Litvak and Ejov \cite{MR2542990} stating that for any Hamiltonian graph $G$, the permutation matrices associated to the Hamiltonian cycles of $G$ are the unique minimizers of $\fF$ on $\cK(G,\upsilon)$.
But Proposition \ref{contdiscr} does not enable us to extend directly Theorem \ref{theo3} to the discrete time setting,
because the diagonal of the generator associated to a Hamiltonian cycle is constant if and only if the underlying invariant probability measure is uniform.
This extension is nevertheless true. To show it, note that the differentiation technique of Section 3 can be adapted to $\cK(G,\pi)$ in a straightforward manner.

\appendix

\section{APPENDIX: Computations on the simplest example of non-Hamiltonian connected graph}

The length 2 segment $S_2\df (\{0,1,2\},\{(0,1),(1,0),(1,2),(2,1)\})$ is the simplest non-Hamiltonian (strongly) connected graph.
We compute here the minimizer of $F$ on $\cL(S_2,\pi)$, for any positive probability measure $\pi$ on $\{0,1,2\}$.
We hope this example will motivate further investigation of the minimizers of $F$ in the challenging non-Hamiltonian framework.
\par\me
To simplify the notation, write $x=\pi(0), y=\pi(1)$ and $z=\pi(2)$, by assumption we have that $x,y,z>0$ and $x+y+z=1$.
Up to exchanging the vertices 0 and 2, we assume that $\vert x-1/2\vert\geq \vert z-1/2\vert$.\par
Any Markov generator $L$ from $\cL(S_2,\pi)$ has the form
\bq
L&\df& \lt(\begin{array}{ccc}
-a&a&0\\
\alpha&-\alpha-\beta& \beta\\
0&b& -b
\end{array}\rt)\eq
where the coefficients $a,\alpha,\beta, b>0$ satisfy,
\bq
xa&=&y\alpha,\\
y\beta&=&z b,\\
xa+y(\alpha+\beta)+zb&=&1.\eq
The first two equalities correspond to the invariance  of $\pi$ for $L$ (here $\pi$ is even reversible for the birth and death generator $L$)
and the third one is the normalization condition, it can be rewritten
\bqn{renorma}
2xa+2zb&=&1\eqn
Denote $\Lambda(L)=\{\lambda_1,\lambda_2\}$, its elements are the non-zero roots in $X$ of the polynomial
$\det(X+L)$. We compute that
\bq
\det(X+L)&=&X(X^2-(a+\alpha+\beta+b)X+ab+a\beta+\alpha b),\eq
so that
\bq
\lambda_1+\lambda_2&=&a+\alpha+\beta+b,\\
\lambda_1\lambda_2&=&ab+a\beta+\alpha b.\eq
From \eqref{spectral1}, we have
\bq
F(L)& =& \frac1{\lambda_1}+\frac1{\lambda_2}\\& =& \frac{\lambda_1+\lambda_2}{\lambda_1\lambda_2}\\
&=&\frac{a+\alpha+\beta+b}{ab+a\beta+\alpha b}\\
&=&\frac{a\lt(1+\frac{x}{y}\rt)+b\lt(1+\frac{z}{y}\rt)}{ab\lt(1+\frac{x}{y}+\frac{z}{y}\rt)}\\
&=&\frac{a(x+y)+b(y+z)}{ab}\\
&=&\frac{x+y}{b}+\frac{y+z}{a}\\
&=&\frac{1-z}{b}+\frac{1-x}{a}.
\eq
Taking into account \eqref{renorma},
the minimizer of $F$ on $\cL(S_2,\pi)$ corresponds to the minimizer of
\bqn{mini}
(0, 1/(2x))\ni a&\mapsto& 2z\frac{1-z}{1-2xa}+\frac{1-x}a.\eqn
We are thus led to the second order equation in $a$:
\bqn{deux}
4x(z(1-z)-x(1-x))a^2+4x(1-x)a-(1-x)&=&0.\eqn
Due to the assumption $\vert x-1/2\vert\geq \vert z-1/2\vert$, the first coefficient is non-negative.
We consider two cases.\par\sm
$\bullet$ If $\vert x-1/2\vert= \vert z-1/2\vert$, then \eqref{deux} degenerates into a first order equation and
$a\df 1/(4x)$ is the minimizer of the mapping \eqref{mini}.
It follows that the minimizer of $F$ on $\cL(S_2,\pi)$ is
\bq
L_\wedge&\df& \setlength{\extrarowheight}{9pt}\lt(\begin{array}{*3{>{\displaystyle}c}}
-\frac1{4x}&\frac1{4x}&0\\
\frac1{4y}&-\frac1{2y}& \frac1{4y}\\
 0&\frac1{4z}& -\frac1{4z}
\end{array}\rt)\eq
and the minimal value $F_\wedge(S_2,\pi)$ of $F$ on $\cL(S_2,\pi)$ is
\bq
F(L_\wedge)\ =\ 4(1-z)z+4(1-x)x\ =\ 8x(1-x).\eq
In particular, for $\pi=\upsilon$, the uniform distribution on $\{0,1,2\}$, we get
\bq
L_\wedge&\df& \frac14 \lt(\begin{array}{ccc}
-3&3&0\\
3&-6& 3\\
0&3& -3
\end{array}\rt)\eq
and  $F_\wedge(S_2,\upsilon)= 16/9$.\par\sm
$\bullet$ If $\vert x-1/2\vert> \vert z-1/2\vert$, then \eqref{deux} admits two
solutions
\bq
a_\pm&\df& \frac{-x(1-x)\pm\sqrt{x(1-x)z(1-z)}}{2x(z(1-z)-x(1-x))},\eq
but only $a_+$ belongs to $(0,1/(2x))$ and is in fact the minimizer of the mapping \eqref{mini}.
This value can be simplified into
\bq
a_+&=&\frac1{2x}\frac{\sqrt{x(1-x)}}{\sqrt{x(1-x)}+\sqrt{z(1-z)}}.\eq
It follows that the minimizer of $F$ on $\cL(S_2,\pi)$ is
\bq
L_\wedge&\df&\setlength{\extrarowheight}{12pt} \frac1{\sqrt{x(1-x)}+\sqrt{z(1-z)}}\lt(\begin{array}{*3{>{\displaystyle}c}}
-\frac{\sqrt{x(1-x)}}{2x}&\frac{\sqrt{x(1-x)}}{2x}&0\\
\frac{\sqrt{x(1-x)}}{2y}&- \frac{\sqrt{x(1-x)}+\sqrt{z(1-z)}}{2y}& \frac{\sqrt{z(1-z)}}{2y}\\
0&\frac{\sqrt{z(1-z)}}{2z}& -\frac{\sqrt{z(1-z)}}{2z}
\end{array}\rt)\\
&=& pL_{(0,1)}+(1-p)L_{(1,2)}
\eq
with the notation introduced in \eqref{LA} and
\bq
p&\df& \frac{\sqrt{x(1-x)}}{\sqrt{x(1-x)}+\sqrt{z(1-z)}}.\eq\par
The minimal value $F_\wedge(S_2,\pi)$ of $F$ on $\cL(S_2,\pi)$ is
\bq
F(L_\wedge)& =&2\lt(\sqrt{x(1-x)}+\sqrt{z(1-z)}\rt)^2.\eq
\par
Letting $\vert x-1/2\vert$ converge to $ \vert z-1/2\vert$, we recover the values of $L_\wedge$ and $F(L_\wedge)$ obtained in the previous case.
\par


 \bibliographystyle{plain}

\vskip2cm
\hskip70mm
\vbox{
\copy5
\vskip5mm
\copy6
}

\end{document}